\documentclass[11pt]{amsart}

\usepackage{amsfonts,amsmath,amssymb,amsthm}

\newtheorem{thm}{Theorem}[section]
\newtheorem{lem}[thm]{Lemma}
\newtheorem{prp}[thm]{Proposition}
\newtheorem{cor}[thm]{Corollary}
\theoremstyle{definition}
\newtheorem{dfn}[thm]{Definition}
\theoremstyle{remark}
\newtheorem{rem}[thm]{Remark}
\newtheorem{ex}[thm]{Example}

\numberwithin{equation}{section}

\DeclareMathOperator{\Ric}{Ric}
\DeclareMathOperator{\grad}{grad}

\DeclareMathOperator{\diver}{div}

\begin{document}
\thispagestyle{empty}
\title{Curvature of Multiply Warped Products}
\author{Fernando Dobarro}
\thanks{Research of F. D. partially supported by funds
of the National Group 'Analisi Reale' of the Italian Ministry of
University and Scientific Research at the University of Trieste.}
\address{Dipartimento di Matematica e Informatica,
Universit\`{a} degli Studi di Trieste,
Via Valerio 12/B I-34127, Trieste \\
Italy} \email {dobarro@dsm.univ.trieste.it}
\author{B\"{u}lent \"{U}nal}
\date{October 26, 2004}
\address{Department of Mathematics, Atilim University,
         Incek, 06836 Ankara, Turkey}
\email{bulentunal@mail.com} 
%\thanks{}

\begin{abstract}\
In this paper, we study Ricci-flat and Einstein Lorentzian
multiply warped products. We also consider the case of having
constant scalar curvatures for this class of warped products.
Finally, after we introduce a new class of space-times called as
generalized Kasner space-times, we apply our results to this kind
of space-times as well as other relativistic space-times, i.e.,
Reissner-Nordstr\"om, Kasner space-times,
Ba\~{n}ados-Teitelboim-Zanelli and de Sitter Black hole solutions.
\end{abstract}

\maketitle

%\tableofcontents

\vspace{1.5 cm}

\centerline{\textsc{Contents}}

\vspace{0.5 cm}

\begin{center}
\begin{minipage}[t]{10.0 cm}
\begin{list}{\item contents}{}
  \item[1.]Introduction
  %\dotfill \pageref{chp1}
  \item[2.]Preliminaries
  %\dotfill \pageref{chp2}
  \item[3.]Special Multiply Warped Products
  %\dotfill \pageref{chp3}
   \begin{itemize}
  \item[3.1.]Einstein Ricci Tensor
  %\dotfill \pageref{chp3a}
  \item[3.2.]Constant Scalar Curvature
  %\dotfill \pageref{chp3b}
   \end{itemize}
  \item[4.]Generalized Kasner Space-time
  %\dotfill \pageref{chp4}
  \item[5.]$4$-Dimensional Space-time Models
  %\dotfill \pageref{chp5}
   \begin{itemize}
     \item[5.1.]Type (I)
     %\dotfill \pageref{chp5a}
     \item[5.2.]Type (II)
     %\dotfill \pageref{chp5b}
     \item[5.3.]Type (III)
     %\dotfill \pageref{chp5c}
\end{itemize}
  \item[6.]BTZ $(2+1)$-Black Hole Solutions
  %\dotfill \pageref{chp6}
  \item[7.]Conclusions
  %\dotfill \pageref{chp7}
  \item[References]
  %\dotfill \pageref{ref}
\end{list}
\end{minipage}
\end{center}

\newpage

\renewcommand{\thepage}{\arabic{page}}
\setcounter{page}{2}

\section{Introduction} \label{chp1}

The concept of warped products was first introduced by Bishop and
O'Neill (see \cite{BO}) to construct examples of Riemannian
manifolds with negative curvature. In Riemannian geometry, warped
product manifolds and their generic forms have been used to
construct new examples with interesting curvature properties since
then (see \cite{EM,BO,CDVV,DD,EJK,GMV,GA,KD,KK,KC,KN,LR}). In
Lorentzian geometry, it was first noticed that some well known
solutions to Einstein's field equations can be expressed in terms
of warped products in \cite{GLG1} and after that Lorentzian warped
products have been used to obtain more solutions to Einstein's
field equations (see \cite{GLG1,GLG,EM,BO,HE,KSHMC,SRG}).
Moreover, geometric properties such as geodesic structure or
curvature of Lorentzian warped products have been studied by many
authors because of their relativistic applications (see
\cite{ARS,AD1,AD,GES,BJ,BJK,BEP,BJPP,CJ,CJ1,SSS,DK,EJK,
EJKS,ES,FS,FS1,HCP,KKK,KO,KO2, MM,PT,RASM,SM0,SM,SM1,DWP,MWP}).

We recall the definition of a warped product of two
pseudo-Riemannian manifolds $(B,g_B)$ and $(F,g_F)$ with a smooth
function $b \colon B \to (0,\infty)$ (see also \cite{GLG,SRG}).
Suppose that $(B,g_B)$ and $(F,g_F)$ are pseudo-Riemannian
manifolds and also suppose that $b \colon B \to (0,\infty)$ is a
smooth function. Then the (singly) warped product, $B \times {}_b
F$ is the product manifold $B \times F$ equipped with the metric
tensor $g=g_B \oplus b^{2}g_F$ defined by
$$g=\pi^{\ast}(g_B) \oplus (b \circ \pi)^2 \sigma^{\ast}(g_F)$$
where $\pi \colon B \times F \to B$ and $\sigma \colon B \times F
\to F$ are the usual projection maps and ${}^\ast$~denotes the
pull-back operator on tensors. Here, $(B,g_F)$ is called as the
base manifold and $(F,g_F)$ is called as the fiber manifold and
also $b$ is called as the warping function.

Generalized Robertson-Walker space-time models (see
\cite{ARS,BJK,FS,RASM,SM,SM1}) and standard static space-time
models (see \cite{AD1,AD,GES,KO,KO2}) that are two well known
solutions to Einstein's field equations can be expressed as
Lorentzian warped products. Clearly, the former is a natural
generalization of Robertson-Walker space-time and the latter is a
generalization of Einstein static universe. One way to generalize
warped products is to consider the case of multi fibers to obtain
more general space-time models (see examples given in {\it Section
2}) and in this case the corresponding product is so called
multiply warped product. In \cite{MWP}, covariant derivative
formulas for multiply warped products are given and the geodesic
equation for these spaces are also considered. The causal
structure, Cauchy surfaces and global hyperbolicity of multiply
Lorentzian warped products are also studied. Moreover, necessary
and sufficient conditions are obtained for null, time-like and
space-like geodesic completeness of Lorentzian multiply products
and also geodesic completeness of Riemanninan multiply warped
products. In \cite{CJ,CJ1}, the author studies manifolds with
$C^{0}$-metrics and properties of Lorentzian multiply warped
products and then he shows a representation of the interior
Schwarzschild space-time as a multiply warped product space-time
with certain warping functions. He also gives the Ricci curvature
in terms of $b_1, b_2$ for a multiply warped product of the form
$M=(0, 2m) \times _{b_1} \mathbb R^1 \times _{b_2} \mathbb S^2$.
In \cite{HCP}, physical properties (2+1) charged
Ba\~{n}ados-Teitelboim-Zanelli (BTZ) black holes and (2+1) charged
de Sitter (dS) black holes are studied by expressing these metric
as multiply warped product space-times, more explicitly, Ricci and
Einstein tensors are obtained inside the event horizons (see also
\cite{BTZ}). In \cite{SM0}, the existence, multiplicity and causal
character of geodesics joining two points of a wide class of
non-static Lorentz manifolds such as intermediate
Reissner-Nordstr\"om or inner Schwarzschild and generalized
Robertson-Walker space-times are studied. In \cite{FS1}, geodesic
connectedness and also causal geodesic connectedness of
multi-warped space-times are studied by using the method of
Brouwer's topological degree for the solution of functional
equations. There are also different types of warped products such
as a kind of warped product with two warping functions acting
symmetrically on the fiber and base manifolds, called as a doubly
warped product (see \cite{DWP}) or another kind of warped product
called as a twisted product when the warping function defined on
the product of the base and fiber manifolds (see \cite{TWP}).
Moreover, Easley studied {\it Local Existence Warped Product
Structures} and also defined and considered another form of a
warped product in his thesis (see \cite{EK}).

In this paper, we answer some questions about the existence of
nontrivial warping functions for which the multiply warped product
is Einstein or has a constant scalar curvature. This problem was
considered especially for Einstein Riemannian warped products with
compact base and some partial answers were also provided (see
\cite{GA,KD,KK, KC}). In \cite{KK}, it is proved that an Einstein
Riemannian warped product with a non-positive scalar curvature and
compact base is just a trivial Riemannian product. Constant scalar
curvature of warped products was studied in \cite{CDM,DD,EJK,EJKS}
when the base is compact and of generalized Robertson-Walker
space-times in \cite{EJK}. Furthermore, partial results for warped
products with non-compact base were obtained in \cite{BD} and
\cite{CB}. The physical motivation of existence of a positive
scalar curvature comes from the positive mass problem. More
explicitly, in general relativity the positive mass problem is
closely related to the existence of a positive scalar curvature
(see \cite{ZX}). As a more general related reference, one can
consider \cite{KJ} to see a survey on scalar curvature of
Riemannian manifolds. The problem of existence of a warping
function which makes the warped product Einstein was already
studied for special cases such as generalized Robertson-Walker
space-times and a table given the different cases of Einstein
generalized Robertson-Walker when the Ricci tensor of the fiber is
Einstein in \cite{ARS} (see also references therein). Einstein
Ricci tensor and constant scalar curvature of standard static
space-times with perfect fluid were already considered in
\cite{KO,MM}. Moreover, in \cite{KO2}, the conformal tensor on
standard static space-times with perfect fluid is studied and it
is shown that a standard static space-time with perfect fluid is
conformally flat if and only if its fiber is Einstein and hence of
constant curvature. In \cite{SSS}, this problem is considered for
arbitrary standard static space-times, more explicitly, an
essential investigation of conditions for the fiber and warping
function for a standard static space-time (not necessarily with
perfect fluid) is carried out so that there exists no nontrivial
function on the fiber guaranteing that the standard static
space-time is Einstein. Duggal studied the scalar curvature of
4-dimensional triple Lorentzian products of the form $L \times B
\times {}_fF$ and obtained explicit solutions for the warping
function $f$ to have a constant scalar curvature for this class of
products (see \cite{DK}). Moreover, in the present paper, we
introduce an original form to generalize Kasner space-times and
then we obtain necessary and sufficient conditions as well as
explicit solutions, for some special cases, for a generalized
Kasner space-time to be Einstein or to have constant scalar
curvature. Besides than the form mentioned here, there are also
other generalizations in the literature (see \cite{IT,KSS}). In
\cite{IT}, an extension for Kasner space-times is introduced in
the view of generalizing 5-dimensional Randall-Sundrum model to
higher dimensions and in \cite{KSS}, another multi-dimensional
generalization of Kasner metric is described and essential
solutions are also obtained for this class of extension. One can
also consider \cite{{Dabrowski99}, {Frolov01},
{Ivashchuk-Singleton04}, {JMS}, {Maceda04}, {Papadopoulos04}, {van
den Hoogen-Horne04}} for recent applications of Kasner metrics and
its generalizations.

We organize the paper as follows. In {\it Section 2}, we give
several basic geometric facts related to the concept of curvatures
(see \cite{PHD,MWP}). Moreover, we recall two well known examples
of relativistic space-times which can be considered as generalized
multiply Robertson-Walker space-times. In {\it Section 3}, we
obtain two results in which, under several assumptions on the
fibers and warping functions, multiply generalized
Robertson-Walker space-times are Einstein or have constant scalar
curvature. In {\it Section 4}, after we introduce generalized
Kasner space-times, we state conditions for this class of
space-times to be Einstein or to have constant scalar curvature.
In {\it Section 5}, we give an explicit classification of
4-dimensional multiply generalized Robertson-Walker space-times
and 4-dimensional generalized Kasner space-times which are
Einstein. In the last section, we focus on BTZ (2+1)- Black Hole
solutions and classify (BTZ) black hole solutions given in {\it
Section 2} by using a more formal approach (see
\cite{BHTZ,BTZ,HCP,MTZ}) and then we also prove necessary and
sufficient conditions for the lapse function of a BTZ (2+1) Black
Hole solution to have a constant scalar curvature or to be
Einstein. Our main results are obtained in {\it Sections 3,4} and
{\it 5,} especially see {\it Theorem \ref{gric-fe}}, {\it
Propositions \ref{ks-1}} and {\it \ref{ks-2}} as well as {\it
Tables 1,2} and {\it 3}.

\section{Preliminaries} \label{chp2}

Throughout this work any manifold $M$ is assumed to be connected,
Hausdorff, paracompact and smooth. Moreover, $I$ denotes for an
open interval in $\mathbb R$ of the form $I=(t_1, t_2)$ where
$-\infty \leq t_1 < t_2 \leq \infty $ and we will furnish $I$ with
a negative metric $-{\rm d}t^2.$ A pseudo-Riemannian manifold
$(M,g)$ is a smooth manifold with a metric tensor $g$ and a
Lorentzian manifold $(M,g)$ is a {\it pseudo-Riemannian} manifold
with signature $(-,+,+,\cdots,+).$ Moreover, we use the definition
and the sign convention for the {\it curvature} as in \cite{GLG}.
For an arbitrary $n$-dimensional pseudo-Riemannian manifold
$(M,g)$ and a smooth function $f \colon M \to \mathbb R,$ we have
that ${\rm H}^f$ and $\Delta(f)$ denote the {\it Hessian (0,2)}
tensor and the Laplace-Beltrami operator of $f,$ respectively
(\cite{SRG}). Here, we use the sign convention for the Laplacian
in \cite{SRG}, i.e., defined by $\Delta={\rm tr}_g({\rm H}),$ (see
page 86 of \cite{SRG}) where ${\rm H}$ is the Hessian form (see
page 86 of \cite{SRG}) and ${\rm tr}_g$ denotes for the trace, or
equivalently, $\Delta={\rm div}({\rm grad}),$ where ${\rm div}$ is
the divergence and ${\rm grad}$ is the gradient (see page 85 of
\cite{SRG}). Furthermore, we will frequently use the notation
$\|\grad f \|^2=g(\grad f ,\grad f)$. When there is a possibility
any misunderstanding, we will explicitly state the manifold or the
metric for which the operator is considered.
%
%\notemarg
%

We begin our discussion by giving the formal definition of a
multiply warped product (see \cite{MWP}).

\begin{dfn}Let $(B,g_B)$ and $(F_i,g_{F_i})$ be {\it pseudo-Riemannian}
manifolds and also let $b_i \colon B \to (0,\infty)$ be smooth
functions for any $i \in \{1,2,\cdots,m\}.$ The {\it multiply
warped product} is the {\it product manifold} $M=B \times F_1
\times F_2 \times \cdots \times F_m$ furnished with the metric
tensor $g=g_B \oplus b_{1}^{2}g_{F_1} \oplus b_{2}^{2}g_{F_2}
\oplus \cdots \oplus b_{m}^{2}g_{F_m}$ defined by
\begin{equation} \label{idwp}
g=\pi^{\ast}(g_B) \oplus (b_1 \circ \pi)^2
\sigma_{1}^{\ast}(g_{F_1}) \oplus \cdots \oplus (b_m \circ
\pi)^2 \sigma_{m}^{\ast}(g_{F_m})
\end{equation}
\end{dfn}

Each function $b_i \colon B \to (0,\infty)$ is called a {\it
warping} functions and also each manifold $(F_i,g_{F_i})$ is
called a fiber manifold for any $i \in \{1,2,\cdots,m\}.$ The
manifold $(B,g_B)$ is the base manifold of the multiply warped
product.
\begin{itemize}
\item If $m=1,$ then we obtain a {\it singly warped product}.
\item If all $b_i \equiv 1,$ then we have a (trivial) {\it product
manifold}.
\item If $(B,g_B)$ and $(F_i,g_{F_i})$ are all {\it Riemannian}
manifolds for any $i \in \{1,2,\cdots,m\},$ then $(M,g)$ is also
a {\it Riemannian} manifold.
\item The multiply warped product $(M,g)$ is a {\it Lorentzian
multiply warped product} if $(F_i,g_{F_i})$ are all {\it
Riemannian} for any $i \in \{1,2,\cdots,m\}$ and either
$(B,g_B)$ is {\it Lorentzian} or else $(B,g_B)$ is a
one-dimensional manifold with a {\it negative definite} metric
$-{\rm d}t^2$.
\item If $B$ is an open connected interval $I$ of the form
$I=(t_1,t_2)$ equipped with the negative definite metric
$g_B=-{\rm d}t^2,$ where $-\infty \leq t_1 < t_2 \leq \infty,$
and $(F_i,g_{F_i})$ is Riemannian for any $i \in
\{1,2,\cdots,m\},$ then the Lorentzian multiply warped product
$(M,g)$ is called a multiply generalized Robertson-Walker
space-time or a multi-warped space-time. In particular, a
multiply generalized Robertson-Walker space-time is called a
generalized Reissner-Nordstr\"om space-time when $m=2.$
\end{itemize}

We will state the covariant derivative formulas for {\it
multiply warped products} (see \cite{CJ,PHD,MWP}).

\begin{prp} \label{gcovd} Let
$M=B \times {}_{b_1}F_1 \times \cdots \times {}_{b_m}F_m$
be a {\it pseudo-Riemannian} multiply warped product with metric
$g=g_B \oplus b_{1}^{2}g_{F_1} \oplus \cdots \oplus b_{m}^{2}g_{F_m}$
also let  $X,Y \in \mathfrak L(B)$ and
$V \in \mathfrak L(F_i),$ $W \in \mathfrak L(F_j).$ Then
\begin{enumerate}
\item ${\displaystyle \nabla_{X} Y= \nabla_{X}^B Y}$
\item ${\displaystyle \nabla_{X} V= \nabla_{V} X=\frac{X(b_i)}{b_i} V }$
\item $ \nabla_{V} W=
  \begin{cases}
{\displaystyle 0}&
\text{if $ i \neq j $ }, \\
{\displaystyle \nabla_{V}^{F_i} W-\frac{g(V,W)}{b_i} {\rm
grad}_{B}b_i } & \text{if $ i=j $ }
\end{cases} \nonumber $
\end{enumerate}
\end{prp}

One can compute the {\it gradient} and the {\it Laplace-Beltrami}
operator on $M$ in terms of the {\it gradient} and the {\it
Laplace-Beltrami} operator on $B$ and $F_i,$ respectively. From
now on, we assume that $\Delta = \Delta_M$ and $\grad = \grad_M$
to simplify the notation.
%
%\notemarg
%

\begin{prp} \label{ggra-lap}
Let $M=B \times {}_{b_1}F_1 \times \cdots
\times {}_{b_m}F_m$ be a {\it pseudo-Riemannian} multiply warped
product with metric $g=g_B \oplus b_{1}^{2}g_{F_1} \oplus \cdots
\oplus b_{m}^{2}g_{F_m}$ and $\phi \colon B \to \mathbb R$
and ${\psi}_i \colon F_i \to \mathbb R$ be smooth functions for any
$i \in \{1,\cdots,m\}.$ Then
\begin{enumerate}
\item $\grad(\phi \circ \pi)=\displaystyle{\grad_{B}\phi}$
\item $\grad({\psi}_i \circ {\sigma}_i)=\displaystyle{
\frac{\grad_{F_i}{\psi}_i}{b_{i}^2}}$
\item $\Delta(\phi \circ \pi)=\displaystyle{\Delta_{B}\phi+
\sum_{i=1}^m s_i \frac{g_{B}(\grad_{B}\phi,\grad_{B}b_i)}{b_i}}$
\item $\Delta({\psi}_i \circ {\sigma}_i) =\displaystyle{
\frac{\Delta_{F_i}{\psi}_i}{b_{i}^2}}$
\end{enumerate}
\end{prp}

Now, we will state {\it Riemannian curvature} and {\it Ricci
curvature} formulas from \cite{PHD}.

\begin{prp}\label{grcur} Let $M=B \times {}_{b_1}F_1 \times \cdots
\times {}_{b_m}F_m$ be a {\it pseudo-Riemannian} multiply warped
product with metric $g=g_B \oplus b_{1}^{2}g_{F_1} \oplus \cdots
\oplus b_{m}^{2}g_{F_m}$ also let $X,Y,Z \in \mathfrak L(B)$ and
$V  \in \mathfrak L(F_i), W  \in \mathfrak L(F_j)$
and $U \in \mathfrak L(F_k).$ Then
\begin{enumerate}
\item $\displaystyle{R(X,Y)Z = R_{B}(X,Y)Z}$
\item $\displaystyle{R(V,X)Y = -\frac{{\rm H}_B^{b_i}(X,Y)}{b_i}V }$
\item $\displaystyle{R(X,V)W=R(V,W)X=R(V,X)W=0}$ if \, $i \neq j.$
\item $\displaystyle{R(X,Y)V=0}$
\item $\displaystyle{R(V,W)X=0}$ if \, $i=j.$
\item $\displaystyle{R(V,W)U=0}$ if \, $i=j$ and $i,j \neq k.$
\item $\displaystyle{R(U,V)W=-g(V,W)\frac{g_{B}(\grad_{B}b_i,
\grad_{B}b_k)}{b_i b_k}U}$ if \, $i=j$ and $i,j \neq k.$
\item $\displaystyle{R(X,V)W=\frac{g(V,W)}{b_i}\nabla^{B}_{X}
(\grad_{B}b_i)}$ if \, $i=j.$
\item $\displaystyle{R(V,W)U= R_{F_i}(V,W)U+
\frac{\| \grad_{B}b_i \|_{B}^{2}} {b_{i}^2} \left(g(V,U)W-g(W,U)V
\right)} \textit{ if } i,j=k.$
\end{enumerate}
\end{prp}

\begin{prp}\label{gricur} Let $M=B \times {}_{b_1}F_1 \times \cdots
\times {}_{b_m}F_m$ be a {\it pseudo-Riemannian} multiply warped
product with metric $g=g_B \oplus b_{1}^{2}g_{F_1} \oplus \cdots
\oplus b_{m}^{2}g_{F_m},$ also let $X,Y,Z \in \mathfrak L(B)$ and
$V  \in \mathfrak L(F_i)$ and $W  \in \mathfrak L(F_j).$ Then
\begin{enumerate}
\item $\displaystyle{\Ric(X,Y)=\Ric_{B}(X,Y)-\sum_{i=1}^m
\frac{s_i}{b_i}{\rm H}_{B}^{b_i}(X,Y)}$
\item $\displaystyle{\Ric(X,V)=0}$
\item $\displaystyle{\Ric(V,W)=0}$ if \, $i \neq j.$
\item $\displaystyle{\Ric(V,W)=\Ric_{F_i}(V,W)-
\Big(\frac{\Delta_{B}b_i}{b_i}+(s_i-1) \frac{\| \grad_{B}b_i
\|_{B}^{2}}{b_{i}^2}}  \\
\displaystyle{+\sum_{
%\substack
{k=1,
%\\
k \neq i }}^m s_k \frac{g_{B}(\grad_{B}b_i,\grad_{B}b_k)}{b_i
b_k}\Bigl)g(V,W)}$ if \, $i=j.$
\end{enumerate}
\end{prp}

Now, we will compute the {\it scalar curvature} of a {\it multiply
warped product}. In order to do that, one can use the following
{\it orthonormal frame} on $M$ constructed as follows:

%
%\notemarg
%

Let $\displaystyle \left\{ \frac{\partial}{\partial x^1},\cdots,
\frac{\partial}{\partial x^r} \right\}$ and $\displaystyle\left\{
\frac{\partial}{\partial y_{i}^1},\cdots, \frac{\partial}
{\partial y_{i}^{s_i}} \right\}$ be orthonormal frames on open
sets $U \subseteq B$ and $V_i \subseteq F_i$, respectively for any
$i \in \{1,\cdots,m\}.$ Then
$$
\left\{ \frac{\partial}{\partial x^1},\cdots, \frac{\partial}
{\partial x^r}, \frac{\partial}{b_1
\partial y_{1}^1},\cdots, \frac{\partial}{b_1 \partial
y_{1}^{s_1}},\cdots, \frac{\partial}{b_m \partial y_{m}^1},\cdots,
\frac{\partial}{b_m \partial y_{m}^{s_m}}\right\}
$$ is an
orthonormal frame on an open set  $W \subseteq B \times F$
contained in $U \times V \subseteq B \times F,$ where $F=F_1
\times \cdots \times F_m.$

\begin{prp}\label{gscamu} Let $M=B \times {}_{b_1}F_1 \times \cdots
\times {}_{b_m}F_m$ be a {\it pseudo-Riemannian} multiply warped
product with metric $g=g_B \oplus b_{1}^{2}g_{F_1} \oplus \cdots
\oplus b_{m}^{2}g_{F_m}.$ Then, $\tau$ admits the following
expressions
\begin{enumerate}
  \item
\begin{eqnarray*}
\tau & = & \tau_{B}-2\sum_{i=1}^m s_i \frac{\Delta_{B}b_i}{b_i}
+\sum_{i=1}^m \frac{\tau_{F_i}}{b_{i}^2}-
\sum_{i=1}^m s_i(s_i-1)\frac{\| \grad_{B}b_i \|_{B}^{2}}{b_{i}^2}\\
& - & \sum_{i=1}^m \sum_{\substack{k=1 , k \neq i }}^m s_k s_i
\frac{g_{B}(\grad_{B}b_i,\grad_{B}b_k)}{b_i b_k},
\end{eqnarray*}
  \item
\begin{eqnarray*}\label{}
\tau  &=&  \tau_{B}-\sum_{i=1}^m s_i \frac{\Delta_{B}b_i}{b_i}
%+\sum_{i=1}^m \frac{\tau_{F_i}}{b_{i}^2}
%
-\diver \sum_{i=1}^m s_i\frac{\grad_{B}b_i}{b_{i}} \\
&-&g_{B}\left(\sum_{i=1}^m s_i
\frac{\grad_{B}b_i}{b_i},\sum_{i=1}^m s_i
\frac{\grad_{B}b_i}{b_i}\right)
+\sum_{i=1}^m \frac{\tau_{F_i}}{b_{i}^2}.
\end{eqnarray*}
\end{enumerate}
\end{prp}

The following formula can be directly obtained from the previous
result and noting that on a multiply generalized Robertson-Walker
space-time
${\rm grad}_{B}b_{i}=-b_{i}^{\prime}$,
$\| \grad_{B}b_i \|_{B}^{2}\\=-(b_{i}^{\prime})^{2}$,
$\displaystyle g_{B}\left(\frac{\partial}{\partial t},
\frac{\partial}{\partial t}\right)=-1$,
%
%$\| \grad_{B}b_i \|_{B}^{2}=-(b_{i}^{\prime})^{2}$,
%
$\displaystyle {\rm H}_{B}^{b_{i}}\left(\frac{\partial}{\partial
t}, \frac{\partial}{\partial t}\right)=b_{i}^{\prime \prime}$ and
$\Delta_{B}b_{i}=-b_{i}^{\prime \prime}$, we denote the usual
derivative on the real interval $I$ by the prime notation (i.e.,
$\prime$) from now on.

%
%\notemarg
%

\begin{cor} \label{gscamu1} Let $M=I \times {}_{b_1}F_1 \times
\cdots \times {}_{b_m}F_m$ be a multiply generalized
Robertson-Walker space-time with the metric $g=-{\rm d}t^2 \oplus
b_{1}^{2}g_{F_1} \oplus \cdots \oplus b_{m}^{2}g_{F_m}.$ Then,
$\tau$ admits the following expressions
\begin{enumerate}
  \item
$$ \tau = 2\sum_{i=1}^m s_i \frac{b^{\prime \prime}_i}{b_i}
+\sum_{i=1}^m \frac{\tau_{F_i}}{b_{i}^2}+ \sum_{i=1}^m
s_i(s_i-1)\frac{(b_i^\prime)^2}{b_{i}^2} + \sum_{i=1}^m
\sum_{\substack{k=1 , k \neq i }}^m s_k s_i \frac{b_i^\prime
b_k^\prime }{b_i b_k}, $$
\item
$$ \tau = \sum_{i=1}^m s_i \frac{b^{\prime \prime}_i}{b_i}
%+\sum_{i=1}^m \frac{\tau_{F_i}}{b_{i}^2}
+ \left(\sum_{i=1}^m s_i\frac{b_i^\prime}{b_{i}}\right)^\prime
+ \left(\sum_{i=1}^m s_i\frac{b_i^\prime}{b_{i}}\right)^2
+\sum_{i=1}^m \frac{\tau_{F_i}}{b_{i}^2}.
 $$
\end{enumerate}
\end{cor}

We now give some physical examples of relativistic space-times and
state some of their geometric properties to stress the physical
motivation and importance of Lorentzian multiply warped products.
The first example is Schwarzschild black hole solution or known as
inner Reissner-Nordstr\"om space-time and the second one is Kasner
space-time. Our last two examples are closely related to each
other, more explicitly, the third example is
Ba\~{n}ados-Teitelboim-Zanelli (BTZ) black hole solution and the
final example is de Sitter (dS) black hole solution.

$\bullet$ {\bf Schwarzschild Space-time}

We will briefly discuss the interior Schwarzschild solution. We
show how the interior solution can be written as a multiply
warped product.

The line element of the {\it Schwarzschild black hole}
space-time model for the region $r<2m$ is given as (see
\cite{HE})
$${\rm d}s^2=-\left({\frac{2m}{r}}-1\right)^{-1}{\rm d}r^2+
\left({\frac{2m}{r}}-1\right) {\rm d}t^2+r^2{\rm d}\Omega^2,$$
where ${\rm d}\Omega^2={\rm d} \theta^2 + \sin^2 \theta {\rm d}
\varphi^2$ on $\mathbb S^2.$

In \cite{CJ}, it is shown that this space-time model can be
expressed as a multiply generalized Robertson-Walker space-time,
i.e.,
$${\rm d}s^2=-{\rm d}t^2+b^2_1(t){\rm d}r^2+b^2_2(t)
{\rm d} \Omega^2,$$ where
$$b_1(t)=\sqrt{{\frac{2m}{F^{-1}(t)}}-1} \quad
\text{and} \quad b_2(t)={F^{-1}(t)} \quad \text{also}$$
$$t=F(r)=2m \arccos \left(\sqrt{\frac{2m-r}{2m}}\:\, \right)-
\sqrt{r(2m-r)} \quad \text{such that}$$
$$\lim_{r \to 2m}F(r)=m \pi \quad \text{and} \quad
\lim_{r \to 0}F(r)=0.$$ Moreover, we also need to impose the
above multiply generalized Robertson-Walker space-time model for
the {\it Schwarzschild black hole} to be Ricci-flat due to the
fact that the {\it Schwarzschild black hole} is Ricci-flat (see
also the review of Miguel S\'anchez in AMS for \cite{CJ}).

$\bullet$ {\bf Kasner Space-time}

We consider the {\it Kasner} space-time as a {\it Lorentzian}
multiply warped product (see \cite{GRA}).

A Lorentzian multiply warped product $(M,g)$ of the form
$M=(0,\infty) \times {}_{t^{p_1}}\mathbb R \times
{}_{t^{p_2}}\mathbb R \times {}_{t^{p_3}}\mathbb R$ with the
metric $g=-{\rm d}t^2 \oplus t^{2p_1}{\rm d}x^2 \oplus
t^{2p_2}{\rm d}y^2 \oplus t^{2p_3}{\rm d}z^2$ is said to be the
{\it Kasner} space-time if $p_1+p_2+p_3=
(p_1)^2+(p_2)^2+(p_3)^2=1$ (see \cite{KA}).

It is known by \cite{HA} that $-1/3 \leq p_1, p_2, p_3 < 1.$ It
is also known that, excluding the case of two $p_i$'s zero, then
one $p_i$ is negative and the other two are positive. Thus we
may assume that $-1/3 \leq p_1 < 0 < p_2 \leq p_3 < 1$ by
excluding the case of two $p_i$'s zero and one $p_i$ equal to 1.
Furthermore, the only solution in which $p_2=p_3$ is given by
$p_1=-1/3$ and $p_2=p_3=2/3.$ Note also that since $-1/3 \leq
p_1, p_2, p_3 < 1,$ we have to assume $B$ to be $(0,\infty).$
Clearly, the {\it Kasner} space-time is {\it globally
hyperbolic} (see \cite{MWP}).

By making use of the results in \cite{MWP}, it can be easily
seen that the {\it Kasner} space-time is future-directed
time-like and future-directed null geodesic complete but it is
past-directed time-like and past-directed null geodesic
incomplete. Moreover, it is also space-like geodesic incomplete.

Notice that the Kasner space-time is Einstein with $\lambda=0$
(i.e., Ricci-flat) (see \cite{KA} and page 135 of \cite{KSHMC})
and hence has constant scalar curvature as zero. This fact can be
proved as a particular consequence of our results in the next
section, namely by using {\it Theorem \ref{gric-fe}}.

$\bullet$ {\bf Static Ba\~{n}ados-Teitelboim-Zanelli (BTZ)
Space-time}

In \cite{HCP}, authors classify (BTZ) black hole solutions in
three different classes as static, rotating and charged. Here, we
will only give a brief description of a static BTZ space-time in
terms of Lorentzian multiply warped products, i.e., multiply
generalized Robertson-Walker space-times (see also
\cite{BHTZ,BTZ,MTZ}). The line element of a static BTZ black hole
solution can be expressed as
$${\rm d}s^2=-N^{-2}{\rm d}r^2+N^2 {\rm d}t^2+r^2 {\rm d}
\Omega^2,$$ where ${\rm d}\Omega^2={\rm d} \theta^2 + \sin^2
\theta {\rm d} \varphi^2$ on $\mathbb S^2.$

The line element of the {\it Static BTZ black hole} space-time
model for the region $r<r_H$ can be obtained by taking
$$N^2=m-\frac{r^2}{l^2}.$$ In this case, the space-time model can be
expressed as a multiply generalized Robertson-Walker space-time,
i.e.,
$${\rm d}s^2=-{\rm d}t^2+b^2_1(t){\rm d}r^2+b^2_2(t)
{\rm d} \Omega^2,$$ where $r_H=l \sqrt m$
$$b_1(t)=\sqrt{m-\frac{(F^{-1}(t))^2}{l^2}} \quad
\text{and} \quad b_2(t)={F^{-1}(t)} \quad \text{also}$$
$$t=F(r)=l \arcsin \left(\frac{r}{r_H}\right)
\quad \text{such that}$$
$$\lim_{r \to r_H}F(r)=\frac{l \pi}{2} \quad \text{and} \quad
\lim_{r \to 0}F(r)=0.$$ Here, note that the constant scalar
curvature $\tau$ of the multiply generalized Robertson-Walker
space-time introduced above is $\tau=-6/l^2$ (see \cite{HCP}) or
apply {\it Corollary \ref{gscamu1}}.

Note that, in \cite{HCP}, they also classify (dS) black hole
solution in three classes as static, rotating and charged,
similar to (BTZ) black hole solutions (see \cite{BHTZ,BTZ,MTZ}).

We now state a couple of results which we will frequently be
applied along this article.

The first one is an easy computation which we will show explicitly
below. Let $(M,g)$ be an $n$-dimensional pseudo-Riemannian
manifold. For any $t \in \mathbb{R}$ and $v \in
C^\infty_{>0}(B)=\{v \in C^\infty(B): v>0\}$,
\begin{equation}
\begin{split}
&\grad_g v^{t} = t  v^{t -1} \grad_g v \\
&\Delta _{g} v^{t} = t [(t-1) v^{t-2} \|\grad_g v\|_g^{2} + v^{t
-1}
\Delta _{g}v]\\
&\frac{\Delta _{g} v^{t}}{v^{t}} = t \left[(t-1) \frac{\|\grad_g
v\|_g^{2}}{v^{2}} + \frac{\Delta_{g}v}{v}\right].
\end{split} \label{m-eq}
\end{equation}

The second one is a lemma that follows (for a proof and some
extensions as well as other useful applications, see Section 2 of
\cite{SBC}).

\begin{lem} \label{m-lem} Let $(M,g)$ be an $n$-dimensional
pseudo-Riemannian manifold. Let $L_{g}$ be a differential operator
on $C^\infty_{>0}(M)$ defined by
\begin{equation} L_{g} v =\displaystyle\sum_{i=1}^k
r_{i}\frac{\Delta_{g}v^{a_{i}}}{v^{a_{i}}}, \label{m-lem1}
\end{equation}

where $r_{i},a_{i} \in \mathbb{R}$ and $\zeta :=\displaystyle
\sum_{i=1}^k r_{i}a_{i}$, $\eta:= \displaystyle \sum_{i=1}^k
r_{i}a_{i}^{2} $. Then,

\begin{itemize}
\item[\bf (i)] \begin{equation} L_{g} v = (\eta -
\zeta)\frac{\|\grad_g v \|_g^{2}}{v^{2}} + \zeta
\frac{\Delta_{g}v}{v}. \label{m-lem2}
\end{equation}

\item[\bf (ii)] If $\zeta \neq 0$ and $\eta \neq 0$, for $\alpha =
\displaystyle \frac{\zeta}{\eta}$ and $\beta = \displaystyle
\frac{\zeta^{2}}{\eta},$ then we have
\begin{equation} L_{g} v= \beta \frac
{\Delta_{g}v^{\frac{1}{\alpha}}}{v^{\frac{1}{\alpha}}}.
\label{m-lem3}
\end{equation}
\end{itemize}
\end{lem}

\section{Special Multiply Warped Products} \label{chp3}

\subsection{Einstein Ricci Tensor} \label{chp3a}

In this section, we state some condition to guarantee that a
{\it multiply generalized Robertson-Walker space-time} is {\it
Ricci-flat} or {\it Einstein}.

Now, we recall some elementary facts about Einstein manifolds
starting from its definition.

Recall that an $n$-dimensional pseudo-Riemannian manifold $(M,g)$
is said to be Einstein if there exists a smooth real-valued
function $\lambda$ on $M$ such that $\Ric=\lambda g,$ and
$\lambda$ is called the Ricci curvature of $(M,g)$ (see also page
7 of \cite{AU}).

\begin{rem} Concerning to this notion, it should be pointed out:
\begin{enumerate}
\item If $(M,g)$ is Einstein and $n \geq 3,$ then $\lambda$ is
constant and $\lambda=\tau/n,$ where $\tau$ is the constant scalar
curvature of $(M,g).$

\item If $(M,g)$ is Einstein and $n=2,$ then $\lambda$ is not
necessarily constant.

\item If $(M,g)$ has constant sectional curvature $k,$ then
$(M,g)$ is Einstein with $\lambda=k(n-1)$ and has constant scalar
curvature $\tau=n(n-1)k.$

\item $(M,g)$ is Einstein with Ricci curvature $\lambda$ and
$n=3,$ then $(M,g)$ is a space of constant (sectional) curvature
${\rm K}=\lambda/2.$

\item If $(M,g)$ is a Lorentzian manifold then $(M,g)$ is Einstein
if and only if $\Ric(v,v)=0,$ for any null vector field $v$ on
$M.$
\end{enumerate}
\end{rem}

By using {\it Proposition \ref{gricur}}, we easily obtain the
{\it Ricci curvature} of {\it Lorentzian multiply warped
products}, $(M,g)$ of the above form.

\begin{prp} \label{gpmriki} Let $M=I \times {}_{b_1}F_1 \times
\cdots \times {}_{b_m}F_m$ be a multiply generalized
Robertson-Walker space-time with the metric $g=-{\rm d}t^2
\oplus b_{1}^{2}g_{F_1} \oplus \cdots \oplus b_{m}^{2}g_{F_m}$
also let $\frac{\partial}{\partial t} \in \mathfrak{X}(I)$ and
$v_i \in \mathfrak{X}(F_i),$ for any $i \in \{1,\cdots,m\}.$ If
$v=\sum_{i=1}^m v_i \in \mathfrak{X}(F),$ then
\begin{eqnarray*}
{\rm Ric}\Big(\frac{\partial}{\partial t}+v,
\frac{\partial}{\partial t}+v\Big) &=& \sum_{i=1}^m \Big({\rm
Ric}_{F_i}(v_i,v_i) +\Bigl( b_i b_{i}^{\prime
\prime}+(s_i-1)(b_{i}^{\prime})^2 \\  &+& b_{i} b_{i}^{\prime}
\sum_{\substack{k=1 , k \neq i }}^m s_k
\frac{b_{k}^{\prime}}{b_{k}} \Bigl) g_{F_i}(v_i,v_i)-
s_i\frac{b_{i}^{\prime \prime}}{b_{i}}\Big)
\end{eqnarray*}
\end{prp}

\begin{proof} By substituting $\displaystyle \overline X=
\frac{\partial}{\partial t} + \sum_{i=1}^m v_i$ and $\displaystyle
\overline Y=\frac{\partial}{\partial t}+ \sum_{i=1}^m v_i$ and by
noting that ${\rm grad}_{B}b_{i}=-b_{i}^{\prime}$, $\displaystyle
g_{B}\left(\frac{\partial}{\partial t}, \frac{\partial}{\partial
t}\right)=-1$, $\displaystyle g_{B}({\rm grad}_{B}b_{i},{\rm
grad}_{B}b_{i})  = -(b_{i}^{\prime})^{2}$, $\displaystyle {\rm
H}_{B}^{b_{i}}\left(\frac{\partial}{\partial t},
\frac{\partial}{\partial t}\right)=b_{i}^{\prime \prime}$,
$\Delta_{B}b_{i}=-b_{i}^{\prime \prime}$
 and $\displaystyle {\rm Ric}_{B}\left(\frac{\partial}{\partial
t}, \frac{\partial}{\partial t}\right)=0$ and by using {\it
Proposition \ref{gricur}}, we obtain the result.
\end{proof}

The following result can be easily proved by substituting $v_j=0$
for any $j \in \{1,\cdots,m\}-\{i\}$ and $v_i \neq 0,$ in {\it
Proposition \ref{gpmriki}} along with the method of separation of
variables.

\begin{thm} \label{gric-fe}  Let $M=I \times {}_{b_1}F_1 \times
\cdots \times {}_{b_m}F_m$ be a multiply generalized
Robertson-Walker space-time with the metric $g=-{\rm d}t^2 \oplus
b_{1}^{2}g_{F_1} \oplus \cdots \oplus b_{m}^{2}g_{F_m}.$ The
space-time $(M,g)$ is Einstein with Ricci curvature $\lambda$ if
and only if the following conditions are satisfied for any $i \in
\{1,\cdots,m\}$
\begin{enumerate}
\item each fiber $(F_i,g_{F_i})$ is Einstein with Ricci curvature
$\lambda_{F_i}$ for any $i \in \{1,\cdots,m\},$ \item
$\displaystyle{\sum_{i=1}^m s_i \frac{b_{i}^{\prime
\prime}}{b_i}=\lambda}$ and \item $\displaystyle{\lambda_{F_i}+b_i
b_{i}^{\prime \prime}+(s_i-1) (b_{i}^{\prime})^2+b_{i}
b_{i}^{\prime} \sum_{\substack{k=1 , k \neq i }}^m s_k
\frac{b_{k}^{\prime}}{b_{k}} = \lambda b_{i}^2}$
\end{enumerate}
\end{thm}

\begin{rem} In {\it Theorem \ref{gric-fe}}, Equation (3) can be
expressed in different forms and here we want to present some of
them. By applying {\it Equation \ref{m-eq}}, we can have
\begin{equation*}\label{}
\frac{\lambda_{F_{i}}}{b_{i}^{2}} + \frac{1}{s_{i}}
\frac{(b_{i}^{s_{i}})^{\prime \prime}}{b_{i}^{s_{i}}} +
\frac{b_{i}^{\prime}}{b_{i}} \sum_{\substack{k=1 , k \neq i }}^{m}
s_{k} \frac{b_{k}^{\prime}}{b_{k}}=\lambda,
%\leqno(3.4-3i)
\leqno(E_{gRW}-i)
\end{equation*}
or equivalently,
\begin{equation*} \label{}
\frac{\lambda_{F_{i}}}{b_{i}^{2}} + \frac{b_{i}^{\prime
\prime}}{b_{i}} -\frac{(b_{i}^{\prime})^{2}}{b_{i}^{2}} +
\frac{b_{i}^{\prime}}{b_{i}}\sum_{k=1}^{m}s_{k}
\frac{b_{k}^{\prime}}{b_{k}} =\lambda.
%\leqno(3.4-3ii)
\leqno(E_{gRW}-ii)
\end{equation*}
\end{rem}

\subsection{Constant Scalar Curvature} \label{chp3b}

It is possible to obtain equivalent expressions for the scalar
curvature in {\it Corollary \ref{gscamu1}}, namely the following
just follows from {\it Equation \ref{m-eq}},
\begin{equation*} \label{}
\tau=\sum_{i=1}^{m}   \left[ s_{i} \frac{b_{i}^{\prime
\prime}}{b_{i}} + \frac{(b_{i}^{s_{i}})^{\prime \prime
}}{b_{i}^{s_{i}}} \right]+ \sum_{i=1}^{m}
\frac{\tau_{F_{i}}}{b_{i}^{2}} + \sum_{i=1}^{m} \sum_{k=1,k \neq
i}^{m} s_{k} s_{i} \frac{b_{i}^{\prime}}{b_{i}}
\frac{b_{k}^{\prime}}{b_{k}}.
%\leqno(3.6-i)
\leqno(sc_{gRW}-i)
\end{equation*}

Since $2 s_{i} \neq 0$ and $ s_{i} + s_{i}^{2} = s_{i}(s_{i} + 1)
\neq 0$, by {\it Lemma \ref{m-lem}}, there results
\begin{equation*} \label{}
\tau = \sum_{i=1}^{m} \frac{4 s_{i}}{s_{i} + 1}
\frac{(b_{i}^{\frac{s_{i} + 1}{2}})^{\prime \prime
}}{b_{i}^{\frac{s_{i} + 1}{2}}} + \sum_{i=1}^{m}
\frac{\tau_{F_{i}}}{b_{i}^{2}} + \sum_{i=1}^{m} \sum_{k=1,k \neq
i}^{m} s_{k} s_{i} \frac{b_{i}^{\prime}}{b_{i}}
\frac{b_{k}^{\prime}}{b_{k}}.
%\leqno(3.6-ii)
\leqno(sc_{gRW}-ii)
\end{equation*}

Thus, defining $\psi_{i} = b_{i}^{\frac{s_{i} + 1}{2}}$, results
\begin{equation*}\label{}
\tau = \sum_{i=1}^{m} \frac{4 s_{i}}{s_{i} + 1}
\frac{\psi_{i}^{\prime \prime}}{\psi_{i}} + \sum_{i=1}^{m}
\frac{\tau_{F_{i}}}{\psi_{i}^{\frac{4}{s_{i}+1}}}+
\sum_{i=1}^{m} \sum_{k=1,k \neq i}^{m} s_{k} s_{i} \frac{(
\psi_{i}^{\frac{2}{s_{i}+1}})^{\prime}}
{\psi_{i}^{\frac{2}{s_{i}+1}}} \frac{(\psi_{k}^{\frac{2}
{s_{k}+1}} )^{\prime}}{\psi_{k}^{\frac{2}{s_{k}+1}}}.
%\leqno(3.6-iii)
\leqno(sc_{gRW}-iii)
\end{equation*}

Note that when $m=1$ this relation is exactly that obtained in
\cite{DD} and \cite{SBC} when the base has dimension $1$.

%Like for Equation $(E_{gRW}-ii)$ of {\it Theorem \ref{gric-fe}},
%it is possible to find another expression of $\tau$.
%\begin{equation*}\label{}
%\tau = \sum_{i=1}^{m} s_{i} \frac{b_{i}^{\prime \prime}}{b_{i}}+
%\sum_{i=1}^{m} s_{i} \left[\frac{b_{i}^{\prime \prime}}{b_{i}}-
%\frac{(b_{i}^{\prime})^{2}}{b_{i}^{2}}\right]+ \sum_{i=1}^{m}
%\frac{\tau_{F_{i}}}{b_{i}^{2}} + \sum_{i=1}^{m} \sum_{k=1}^{m}
%s_{k} s_{i} \frac{b_{i}^{\prime}}{b_{i}}
%\frac{b_{k}^{\prime}}{b_{k}}.
%\leqno(3.6-iv)
%\leqno(sc_{gRW}-iv)
%\end{equation*}

The following result just follows from the method of separation
of variables and the fact that each $\tau_{F_i} \colon F_i \to
\mathbb R$ is a function defined on $F_i,$ for any $i \in
\{1,\cdots,m\}.$

\begin{prp}\label{gscamu2} Let $M=I \times {}_{b_1}F_1 \times
\cdots \times {}_{b_m}F_m$ be a multiply generalized
Robertson-Walker space-time with the metric $g=-{\rm d}t^2
\oplus b_{1}^{2}g_{F_1} \oplus \cdots \oplus b_{m}^{2}g_{F_m}.$
If the space-time $(M,g)$ has constant scalar curvature $\tau,$
then each fiber $(F_i,g_{F_i})$ has constant scalar curvature
$\tau_{F_i},$ for any $i \in \{1,\cdots,m\}.$
\end{prp}

As one can notice from the above formula, it is extremely hard to
determine general solutions for warping functions which produce an
Einstein, or with constant scalar curvature multiply generalized
Robertson-Walker space-time. Note that non-linear second order
differential equations need to be solved according {\it Theorem
\ref{gric-fe}}. Further note that there is only one differential
equation and $m$ different warping functions in {\it Corollary
\ref{gscamu1}}. Therefore instead of giving a general answer to
the existence of warping functions to get an Einstein, or with
constant scalar curvature, space-time, we simplify this problem
and consider some specific cases in mentioned {\it Sections 4} and
{\it 5}.

\section{Generalized Kasner Space-time} \label{chp4}

In this section we give an extension of Kasner space-times and
consider their scalar and Ricci curvatures.

\begin{dfn} \label{gkst} A generalized Kasner space-time $(M,g)$
is a Lorentzian multiply warped product of the form $M=I \times
{}_{\varphi^{p_1}}F_1 \times \cdots \times
{}_{\varphi^{p_m}}F_m$ with the metric $g=-{\rm d}t^2 \oplus
\varphi^{2p_1}g_{F_1} \oplus \cdots \oplus
\varphi^{2p_m}g_{F_m},$ where $\varphi \colon I \to (0,\infty)$
is smooth and $p_i \in \mathbb R,$ for any $i \in
\{1,\cdots,m\}$ and also $I=(t_1, t_2)$ with $-\infty \leq t_1 <
t_2 \leq \infty.$
\end{dfn}

Notice that a Kasner space-time can be obtained out of a form
defined above by taking $\varphi=Id_{(0,\infty)}$ with $m=3$ and
$I=(0,\infty),$ where $Id_{(0,\infty)}$ denotes for the identity
function on $(0,\infty)$ (see \cite{HA}).

From now on, for an arbitrary generalized Kasner space-time of the
form in {\it Definition \ref{gkst}}, we introduce the following
parameters
\begin{equation*} \label{}
\zeta:=\sum_{l=1}^{m} s_{l} p_{l} \quad \textrm{and} \quad
\eta:=\sum_{l=1}^{m} s_{l} p_{l}^{2}.\leqno(\zeta;\eta)
\end{equation*}

\begin{rem} \label{rem-i} Note that $\zeta \neq 0$ implies
$\eta \neq 0$ and in this case, defining $\displaystyle S=
\sum_{l=1}^{m} s_{l}$, results $\displaystyle
\frac{\eta}{\zeta^{2}}\ge\frac{1}{S}$. The latter is for example
consequence of the H\"{o}lder inequality (compare with page 186
of \cite{F99}).
\end{rem}

By applying {\it Theorem \ref{gric-fe}}, we can easily state the
following result and later we will examine the solvability of
the differential equations therein.

\begin{prp} \label{ks-1} Let $M=I \times {}_{\varphi_1} F_1
\times \cdots \times {}_{\varphi_m}F_m$ be a generalized Kasner
space-time with the metric $g=-{\rm d}t^2 \oplus
\varphi^{2p_1}g_{F_1} \oplus \cdots \oplus \varphi^{2p_m}g_{F_m}.$
Then the space-time $(M,g)$ is Einstein with Ricci curvature
$\lambda$ if and only if
\begin{enumerate}
\item each fiber $(F_{i},g_{F_{i}})$ is Einstein with Ricci
curvature $\lambda_{F_{i}}$ for any $i \in \{1, \cdots,m\}$,

\item $\displaystyle{\lambda = \sum_{l=1}^{m} s_{l}
\frac{(\varphi^{p_{l}})^{\prime \prime}}{\varphi^{p_{l}}}=(\eta -
\zeta)\frac{(\varphi^{\prime})^{2}}{\varphi^{2}}+ \zeta
\frac{\varphi^{\prime \prime}}{\varphi}} \textit{ and } $

\item $\displaystyle{\frac{\lambda_{F_{i}}}{\varphi^{2p_{i}}}
+p_{i}\left[(\zeta-1)\frac{(\varphi^{\prime})^{2} }{\varphi^{2}}+
\frac{\varphi^{\prime \prime}}{\varphi}\right]=\lambda}$.
\end{enumerate}
\end{prp}

\begin{rem} \label{rem3}
Moreover, if in {\it Proposition \ref{ks-1}} we assume that $\zeta
\neq 0$ also, then by {\it Remark \ref{rem-i}} is $\eta \neq 0$.
Hence, $(3)$ is equivalent to
\begin{equation*}\label{}
\frac{\lambda_{F_{i}}}{\varphi^{2p_{i}}} + \frac{p_{i}}{\zeta}
\frac{(\varphi^{\zeta})^{\prime \prime}}{\varphi^{\zeta}}
=\lambda,
%\leqno(4.3-2i)
\leqno(E_{gK}^{(3)}-i)
\end{equation*}
and (2) is equivalent to
\begin{equation*}\label{}
\lambda=\frac{\zeta^{2}}{\eta}
\frac{(\varphi^{\frac{\eta}{\zeta}})^{\prime \prime}}
{\varphi^{\frac{\eta}{\zeta}}}.
%\leqno(4.3-3i)
\leqno(E_{gK}^{(2)}-i)
\end{equation*}
\end{rem}

\begin{proof} (of \textit{Proposition \ref{ks-1}} and
\textit{Remark \ref{rem3}})
In order to prove $(3)$, note that Equation $(E_{gRW}-i)$
says \begin{equation*}\label{}
\frac{\lambda_{F_{i}}}{\varphi^{2p_{i}}}+\frac{1}{s_{i}}
\frac{(\varphi^{p_{i}s_{i}})^{\prime \prime}}
{\varphi^{p_{i}s_{i}}}+ \frac{(\varphi^{p_{i}})^{\prime}}
{\varphi^{p_{i}}} \sum_{k=1,k \neq i}^{m} s_{k}
\frac{(\varphi^{p_{k}})^{\prime}}{\varphi^{p_{k}}}=\lambda.
\end{equation*}
Hence, by {\it Equation \eqref{m-eq}},
\begin{equation*}\label{}
\frac{\lambda_{F_{i}}}{\varphi^{2p_{i}}}+p_{i}(p_{i}s_{i}-
1)\frac{(\varphi^{\prime})^{2} }{\varphi^{2}}+
p_{i}\frac{\varphi^{\prime \prime}}{\varphi}+p_{i}
\frac{\varphi^{\prime}}{\varphi} \sum_{k=1,k \neq i}^{m} s_{k}
p_{k} \frac{\varphi^{\prime}}{\varphi}=\lambda,
\end{equation*}
and from here
\begin{equation*}\label{}
\frac{\lambda_{F_{i}}}{\varphi^{2p_{i}}}+ p_{i}
\left[(p_{i}s_{i}-1)+\sum_{k=1,k \neq i}^{m} s_{k} p_{k}
\right]\frac{(\varphi^{\prime})^{2} }{\varphi^{2}}+
p_{i}\frac{\varphi^{\prime \prime}}{\varphi}=\lambda.
\end{equation*}
So,
\begin{equation*}\label{}
\frac{\lambda_{F_{i}}}{\varphi^{2p_{i}}} + p_{i}\left[\left(-1+
\sum_{k=1}^{m} s_{k} p_{k} \right)\frac{(\varphi^{\prime})^{2}
}{\varphi^{2}}+\frac{\varphi^{\prime
\prime}}{\varphi}\right]=\lambda,
\end{equation*}
and by the definition of $\zeta$
\begin{equation*}\label{}
\frac{\lambda_{F_{i}}}{\varphi^{2p_{i}}} + p_{i}\left[(\zeta-1
)\frac{(\varphi^{\prime})^{2} }{\varphi^{2}}+
\frac{\varphi^{\prime \prime}}{\varphi}\right]=\lambda.
\end{equation*}

If furthermore $\zeta \neq 0$, applying again {\it Equation
\ref{m-eq}}, results $(E_{gK}^{(3)}-i)$.

On the other hand, from (2) of {\it Theorem \ref{gric-fe}}
\begin{equation*} \lambda = \sum_{l=1}^{m} s_{l}
\frac{(\varphi^{p_{l}})^{\prime \prime}}{\varphi^{p_{l}}},
\end{equation*}
and by {\it Lemma \ref{m-lem} (a)} ,
\begin{equation*}
\lambda=(\eta-\zeta)\frac{(\varphi^{\prime})^{2}}{\varphi^{2}} +
\zeta \frac{\varphi^{\prime \prime}}{\varphi}.
\end{equation*}
Hence, if $\zeta \neq 0$ and as consequence $\eta \neq 0$,
applying {\it Lemma \ref{m-lem} (b)}, results $(E_{gK}^{(2)}-i)$.
\end{proof}

Note that, from now on and also including the previous result,
when we apply {\it Lemma \ref{m-lem}}, we denote the usual
derivative in equations by means of the prime notation.

\begin{rem} \label{rem4} Note that the conditions $\zeta \neq 0$
and $\eta \neq 0$ agree with the conditions usually imposed in
the classical Kasner space-times, namely $p_{1}+p_{2}+p_{3}=1$
and $p_{1}^{2}+p_{2}^{2}+p_{3}^{2}=1$ (see \cite{KA}). It is
easy to show that the unique possibility to construct an
Einstein classical Kasner manifold or a constant scalar
curvature classical Kasner manifold with $p_{1}+p_{2}+p_{3}=0$
is $p_{1}=p_{2}=p_{3}=0$, so that we have just a usual product.
Indeed, considering $\varphi(t)=t$, it is possible to apply {\it
Proposition \ref{ks-1}} and later {\it Proposition \ref{ks-2}},
respectively.
\end{rem}

\begin{cor} \label{cor1} Under the hypothesis of {\it
Proposition \ref{ks-1}}, along with $\zeta \neq 0$ and $\eta \neq
0$. Assume also that for all $i$, \, $\zeta - p_{i} \neq 0$ and
$\eta - p_{i} \zeta \neq 0.$ Then, $M$ is Einstein if and only if
for any $i \in \{1, \cdots,m\}$, $(F_{i},g_{F_{i}})$ is Einstein
Ricci curvature $\lambda_{F_{i}}$ and
\begin{equation}\label{eq:gKasner-Einstein}
\frac{(\zeta-p_{i})^{2}}{\eta-p_{i}\zeta} \frac{\psi^{\prime
\prime}}{\psi}= \frac{\lambda_{F_{i}}}{\psi^{\frac{\zeta -p_{i}
}{\eta-p_{i}\zeta}2p_{i}}},
%\leqno(4.5)
\end{equation}
where $\displaystyle{0<\psi:=\varphi^{\frac{\eta -
p_{i}\zeta}{\eta-p_{i}}}}$.
\end{cor}

\begin{proof}
Indeed, from equations $(E_{gK}^{(3)}-i)$ and $(E_{gK}^{(2)}-i)$,
\begin{equation*}\label{}
\frac{\lambda_{F_{i}}}{\varphi^{2p_{i}}} =
\frac{\zeta^{2}}{\eta}
\frac{(\varphi^{\frac{\eta}{\zeta}})^{\prime \prime
}}{\varphi^{\frac{\eta}{\zeta}}} -\frac{p_{i}}{\zeta}
\frac{(\varphi^{\zeta})^{\prime \prime}}{\varphi^{\zeta}}.
\end{equation*}
Thus, since for all $i$, \, $\zeta - p_{i} \neq 0$ and $\eta -
p_{i} \zeta \neq 0,$ then applying {\it Lemma \ref{m-lem}}, the
result just follows.
\end{proof}

\begin{ex} \label{rem6} Under the conditions of the classical Kasner
metrics, $m=3$, $p_{1}+p_{2}+p_{3}=1$ and
$p_{1}^{2}+p_{2}^{2}+p_{3}^{2}=1,$ we have $\lambda_{F_{i}}=0$,
$\zeta=1$ and $\eta=1$. Hence the hypothesis $\zeta - p_{i} \neq
0$ and $\eta - p_{i} \zeta \neq 0,$ for all $i,$ implies that
$p_{i}\neq 1$ for all $i$. In this case, Equation
\eqref{eq:gKasner-Einstein} is equivalent to $0<\psi=\varphi$ and
$\psi^{\prime \prime}=0,$ i.e., $0<\varphi(t)=at+b$ with $a,b \geq
0$ and $a^2+b^2 \neq 0.$ Hence, from Equation $(E_{gK}^{(2)}-i)$,
$(0,+\infty)
\times_{\varphi^{p_{1}}}\mathbb{R}\times_{\varphi^{p_{2}}}
\mathbb{R}\times_{\varphi^{p_{3}}}\mathbb{R}$ is Ricci flat
space-time.
\end{ex}

\begin{cor} \label{cor2} Let us assume the hypothesis of {\it Corollary
\ref{cor1}} and that for all $i$, $(F_{i},g_{F_{i}})$ is Ricci
flat. Then, $M$ is Einstein if and only if $\psi^{\prime
\prime}=0$ with $$ 0<\psi:=\varphi^{\frac{\eta -
p_{i}\zeta}{\eta - p_{i}}}, \quad \text{for all} \quad i.$$
\end{cor}

\begin{proof}
It is an immediate consequence of {\it Corollary \ref{cor1}}
\end{proof}

\begin{cor} \label{cor3} Assume that $(F_{i},g_{F_{i}})$ is Ricci
flat for all $i$. Let also $\overline{\zeta},\overline{\eta} \in
\mathbb{R}\setminus\{0\}$ such that $\overline{\zeta}^{2}=
\overline{\eta}$ and $\psi(t)=a t + b \textrm{ with } a,b \geq 0
\textrm{ and } a^2+b^2 > 0.$ If $\zeta=\overline{\zeta}$,
$\eta=\overline{\eta} $, $\zeta-p_{i}\neq 0$ and
$\eta-p_{i_{}}\zeta \neq 0$ for all $i$, then $M = (0,\infty)
\times F_{1} \times_{\varphi^{p_{1}}} \cdots
\times_{\varphi^{p_{m}}} F_{m}$ is a Ricci flat space-time, where
$\varphi=\psi^{\frac{1}{\zeta}}$.
\end{cor}

\begin{proof}
It is sufficient to apply {\it Corollary \ref{cor2}} and {\it
Proposition \ref{ks-1}}.
\end{proof}

\begin{rem} Note that {\it Corollary \ref{cor3}} contains the classical
Kasner metrics except the case in which at least one $p_{i}=1$
(really at most one could be $1$ because
$\eta=p_{1}^{2}+p_{2}^{2}+p_{3}^{2}=1$).
\end{rem}

The following just follows from {\it Corollary \ref{gscamu1}} and
again we discuss the existence of a solution for the differential
equation below.

\begin{prp} \label{ks-2} Let $M=I \times {}_{\varphi_1}
F_1 \times \cdots \times {}_{\varphi_m}F_m$ be a generalized
Kasner space-time with the metric $g=-{\rm d}t^2 \oplus
\varphi^{2p_1}g_{F_1} \oplus \cdots \oplus \varphi^{2p_m}g_{F_m}.$
Then the space-time $(M,g)$ has constant scalar curvature $\tau$
if and only if
\begin{enumerate}
\item each fiber $(F_i,g_{F_i})$ has constant scalar curvature
$\tau_{F_i}$ for any $i \in \{1,\cdots,m\},$ and
\item $\displaystyle{\tau=2\zeta \frac{\varphi^{\prime \prime}}
{\varphi}+ [(\zeta-2)\zeta + \eta] \frac{(\varphi^{\prime})^{2}}
{\varphi^{2}} + \sum_{i=1}^{m}
\frac{\tau_{F_{i}}}{\varphi^{2p_{i}}}},$
\end{enumerate}
\end{prp}

\begin{rem} \label{rem:scgK}If $\zeta \neq 0$, then (2) in
{\it Proposition \ref{ks-2}}
is equivalent to $$\tau = \frac{4\zeta^{2}}{\zeta^{2}+\eta}
\frac{(\varphi^{\frac{\zeta^{2}+\eta}{2\zeta}})^{\prime \prime}}
{\varphi^{\frac{\zeta^{2}+\eta}{2\zeta}}}+ \sum_{i=1}^{m}
\frac{\tau_{F_{i}}}{\varphi^{2p_{i}}}.$$
\end{rem}

\begin{proof} (of \textit{Proposition \ref{ks-2}} and \textit{Remark
\ref{rem:scgK}}) For each $i \in \{1, \cdots,m\}$, let
$\displaystyle \gamma_{i}=p_{i}\frac{s_{i}+1}{2}$ and
$\psi_{i}=\varphi^{\gamma_{i}}$, then by $(sc_{gRW}-iii)$ and {\it
Equation \ref{m-eq}} there results
\begin{eqnarray*}
\tau & = & \sum_{i=1}^{m}\frac{4s_{i}}{s_{i}+1} \gamma_{i}
\left[(\gamma_{i}-1)\frac{(\varphi^{\prime})^{2} }{\varphi^{2}}
+\frac{\varphi^{"}}{\varphi}\right]+ \sum_{i=1}^{m}
\frac{\tau_{F_{i}}}{\varphi^{\frac{4}{s_{i}+1}\gamma_{i}}} \\
& + & \sum_{i=1}^{m} \sum_{k=1,k \neq i}^{m} s_{k} s_{i}
\frac{2\gamma_{i}}{s_{i}+1}
 \frac{2\gamma_{k}}{s_{k}+1}
\displaystyle \frac{(\varphi^{\prime})^{2} }{\varphi^{2}}
\end{eqnarray*}
Then we have
\begin{eqnarray*}
\tau & = & \sum_{i=1}^{m}2s_{i}p_{i} \left[
\left(p_{i}\frac{s_{i}+1}{2}-1\right)\frac{
(\varphi^{\prime})^{2} }{\varphi^{2}}+ \frac{\varphi^{\prime
\prime}}{\varphi}\right]+ \sum_{i=1}^{m}
\frac{\tau_{F_{i}}}{\varphi^{2p_{i}}} \\
& + & \sum_{i=1}^{m} \sum_{k=1,k \neq i}^{m} s_{k} s_{i} p_{i}
p_{k} \frac{(\varphi^{\prime})^{2} }{\varphi^{2}} \\
& = & 2\zeta \frac{\varphi^{\prime \prime}}{\varphi}+
\sum_{i=1}^{m} s_{i}p_{i}\left[2 \left(p_{i}\frac{s_{i}+1}{2}-1
\right)+\sum_{k=1,k \neq i}^{m} s_{k} p_{k}\right]
\frac{(\varphi^{\prime})^{2} }{\varphi^{2}} \\
& + & \sum_{i=1}^{m} \frac{\tau_{F_{i}}}{\varphi^{2p_{i}}} \\
& =& 2\zeta \frac{\varphi^{\prime \prime}}{\varphi}+
\sum_{i=1}^{m} s_{i}p_{i}\left[(\zeta-2) + p_{i}\right]
\displaystyle \frac{(\varphi^{\prime})^{2} }{\varphi^{2}} +
\sum_{i=1}^{m} \frac{\tau_{F_{i}}}{\varphi^{2p_{i}}} \\
& = & 2\zeta \frac{\varphi^{\prime \prime}}{\varphi}+ [(\zeta -
2)\zeta + \eta] \displaystyle \frac{(\varphi^{\prime})^{2}
}{\varphi^{2}} + \sum_{i=1}^{m}
\frac{\tau_{F_{i}}}{\varphi^{2p_{i}}}.
\end{eqnarray*}

Since $ (\zeta - 2)\zeta + \eta + 1 = (\zeta - 1)^{2} + \eta = 0
\textrm{ if and only if } p_{i}=0 \textrm{ for all } i \in \{1,
\cdots,m\}, $ if at least one $p_{i}\neq 0$ there results by {\it
Equation \ref{m-eq}} $$ \tau=(2\zeta-1) \frac{\varphi^{\prime
\prime}}{\varphi}+ \displaystyle \frac{1}{(\zeta - 1)^{2} + \eta}
\displaystyle \frac{(\varphi^{(\zeta-1)^{2} + \eta})^{\prime
\prime} }{\varphi^{(\zeta - 1)^{2} + \eta}} + \sum_{i=1}^{m}
\frac{\tau_{F_{i}}}{\varphi^{2p_{i}}}.
$$ Hence, if $\zeta \neq 0$, applying {\it Lemma
\ref{m-lem}}, $$ \tau = \frac{4\zeta^{2}}{\zeta^{2}+\eta}
\frac{(\varphi^{\frac{\zeta^{2}+\eta}{2\zeta}})^{\prime \prime}}
{\varphi^{\frac{\zeta^{2}+\eta}{2\zeta}}}+ \sum_{i=1}^{m}
\frac{\tau_{F_{i}}}{\varphi^{2p_{i}}}. $$
\end{proof}

\begin{cor} \label{cor13} Under the hypothesis of
{\it Proposition \ref{ks-2}} and $\zeta \neq 0$. Then, by changing
variables as $u=\varphi^{\frac{\zeta^{2}+\eta} {2\zeta}}$, we
conclude that the space-time $M$ has constant scalar curvature
$\tau$ if and only if
\begin{equation*}\label{}
\tau =  \frac{4\zeta^{2}}{\zeta^{2}+\eta} \frac{u^{\prime
\prime}}{u}+ \sum_{i=1}^{m}
\frac{\tau_{F_{i}}}{u^{\frac{4\zeta}{\zeta^{2}+\eta}p_{i}}}
\end{equation*}
or equivalently
\begin{equation*}\label{}
-\frac{4}{1+\displaystyle\frac{\eta}{\zeta^{2}}} u^{\prime
\prime}= -\tau u+ \sum_{i=1}^{m}
\tau_{F_{i}}u^{1-\frac{4}{1+\frac{\eta}{\zeta^{2}}}\frac{p_{i}}{\zeta}}.
\end{equation*}
\end{cor}

\begin{rem} If $\zeta \neq 0$ and there is only one fiber, i.e.,
in a standard warped product, the equation in the previous
corollary corresponds to those obtained in \cite{DD,SBC}.
\end{rem}

\begin{ex} \label{ex1} Let us assume that $\zeta \neq 0$ and
each $F_{i}$ is scalar flat, namely $\tau_{F_{i}}=0.$ Hence,
equation in the previous corollary is written as $$
-\frac{4\zeta^{2}}{\zeta^{2}+\eta} u^{\prime \prime}=-\tau u.$$
Thus all the solutions have the form
\begin{equation*}
u(t)= \begin{cases} \mathcal{A}
e^{i\sqrt{-\tau\frac{\zeta^{2}+\eta}{4\zeta^{2}}}\, t} +
\mathcal{B} e^{-i\sqrt{-\tau\frac{\zeta^{2}+\eta}{4\zeta^{2}}}\,
t} & \text{ if }\tau < 0, \\
\mathcal{A}t + \mathcal{B} & \text{ if }\tau =0,\\
\mathcal{A} e^{\sqrt{\tau\frac{\zeta^{2}+\eta}{4\zeta^{2}}}\, t}
+ \mathcal{B} e^{-\sqrt{\tau\frac{\zeta^{2}+\eta}{4\zeta^{2}}}\,
t} & \text{ if }\tau > 0,
\end{cases}
\end{equation*}
with constants $\mathcal{A}$ and $\mathcal{B}$ such that $u>0$.

If $\zeta = 0$, by {\it Proposition \ref{ks-2}}, we look for
positive solutions of the equation
\begin{equation*}
\begin{cases}
\tau = \eta \displaystyle\frac{(\varphi^\prime)^{2}}{\varphi^{2}}, \\
\varphi > 0.
\end{cases}
\end{equation*}
Since $\eta > 0$, the latter is equivalent to
\begin{equation*}
\begin{cases}
\left(\varphi\displaystyle\frac{\sqrt{\tau}}{\sqrt{\eta}}+
\varphi^{\prime}\right)\left(\varphi\displaystyle\frac{\sqrt{\tau}}
{\sqrt{\eta}}-\varphi^{\prime}\right)=0, \\
\varphi > 0.
\end{cases}
\end{equation*}
Solutions of the equation above are given as,
\begin{equation*}
\varphi(t)=C e^{\pm\frac{\sqrt{\tau}}{\sqrt{\eta}} t},
\end{equation*}
where $C$ is a positive constant.
\end{ex}

Note that this example include the situation of the classical
Kasner space-times in the framework of scalar curvature. Compare
with the results about Einstein classical Kasner metrics in {\it
Remark \ref{rem4}} and {\it Example \ref{rem6}}.

\section{4-Dimensional Space-time Models} \label{chp5}

We first give a classification of 4-dimensional warped product
space-time models and then consider Ricci tensors and scalar
curvatures of them.

\begin{dfn} Let $M=I \times {}_{b_1}F_1 \times \cdots \times
{}_{b_m}F_m$ be a multiply generalized Robertson-Walker space-time
with metric $g=-{\rm d}t^2 \oplus b_{1}^{2}g_{F_1} \oplus \cdots
\oplus b_{m}^{2}g_{F_m}.$
\begin{itemize}
\item $(M,g)$ is said to be of {\it Type (I)}\, if $m=1$ and
${\rm dim}(F)=3.$
\item $(M,g)$ is said to be of {\it Type (II)}\, if $m=2$ and
${\rm dim}(F_1)=1$ and ${\rm dim}(F_2)=2.$
\item $(M,g)$ is said to be of {\it Type (III)}\, if $m=3$ and
${\rm dim}(F_1)=1,$ \, ${\rm dim}(F_2)=1$ and ${\rm dim}(F_3)=1$
\end{itemize}
\end{dfn}

Note that Type (I) contains the Robertson-Walker space-time. The
Schwarzschild black hole solution can be considered as an example
of Type (II). Type (III) includes the Kasner space-time.

\subsection{Type (I)} \label{chp5a}

Let $M=I \times_{b}F$ be a Type (I) warped product space-time with
metric $g=-{\rm d}t^2 \oplus b^{2}g_{F}.$ Then the scalar
curvature $\tau$ of $(M,g)$ is given as

$$ \tau = \frac{\tau_{F}}{b^2}+6\left(\frac{b^{\prime \prime}}{b}
+\frac{(b^\prime)^2}{b^2}\right).$$

The problem of constant scalar curvatures of this type of warped
products, known as generalized Robertson-Walker space-times is
studied in \cite{EJK}, indeed, explicit solutions to warping
function are obtained to have a constant scalar curvature.

If $v$ is a vector field on $F$ and $\displaystyle{\overline{x}
=\frac{\partial}{\partial t}+v},$ then $$\Ric(\overline x,
\overline x)= \Ric_F(v,v)+\bigl(bb^{\prime \prime}+2(b^\prime)^2
\bigl)g_F(v,v)-3 \frac{b^{\prime \prime}}{b}.$$

In \cite{ARS}, explicit solutions are also obtained for the
warping function to make the space-time as Einstein when the
fiber is also Einstein.

\subsection{Type (II)} \label{chp5b}

Let $M=I \times {}_{b_1}F_1 \times {}_{b_2}F_2$ be a Type (II)
warped product space-time with metric $g=-{\rm d}t^2 \oplus
b_{1}^{2}g_{F_1} \oplus b_{2}^{2}g_{F_2}.$ Then the scalar
curvature $\tau$ of $(M,g)$ is given as

$$\tau=\frac{\tau_{F_2}}{b_2^2}+
2\frac{b_1^{\prime \prime}}{b_1}+4\frac{b_2^{\prime
\prime}}{b_2}+2\left(\frac{b_2^\prime}{b_2}\right)^2+
4\frac{b_1^\prime b_2^\prime}{b_1 b_2}.$$

Note that $\tau_{F_1}=0,$ since ${\rm dim}(F_1)=1.$

If $v_i$ is a vector field on $F_i,$ for any $i \in \{1,2\}$ and
$\displaystyle{\overline x=\frac{\partial}{\partial
t}+v_1+v_2},$ then \begin{eqnarray*} \Ric(\overline x,\overline
x) & = & \Ric_{F_2}(v_2,v_2) - \frac{b_1^{\prime
\prime}}{b_1}-\frac{b_2^{\prime \prime}}{b_2} \\
& + & \left( b_1 b_1^{\prime \prime} +
2\frac{b_1 b_1^\prime b_2^\prime}{b_2} \right) g_{F_1}(v_1,v_1) \\
& + & \left( b_2 b_2^{\prime \prime}+(b_2^\prime)^2 + \frac{b_2
b_2^\prime b_1^\prime}{b_1} \right) g_{F_2}(v_2,v_2)
\end{eqnarray*}

Note that $\Ric_{F_1} \equiv 0,$ since ${\rm dim}(F_1)=1.$

$\bullet$ {\bf Classification of Einstein Type (II) generalized
Kasner space-times:}

Let $M = I \times_{\varphi^{p_{1}}}F_{1}\times_{\varphi^{p_{2}}}
F_{2}$ be an Einstein Type (II) generalized Kasner space-time.
Then the parameters introduced before {\it Proposition \ref{ks-1}}
are given by $\zeta = p_{1}+2p_{2}$, $\eta =
p_{1}^{2}+2p_{2}^{2}$. Hence the latter arises

\begin{equation*}
  \begin{cases}
\displaystyle (\eta -
\zeta)\frac{(\varphi^{\prime})^{2}}{\varphi^{2}} + \zeta
\frac{\varphi^{\prime \prime}}{\varphi} =
\lambda  \hfill \\
\displaystyle
p_{1}\left[(\zeta-1)\frac{(\varphi^{\prime})^{2}}{\varphi^{2}} +
\frac{\varphi^{\prime \prime}}{\varphi}
\right]=\lambda\\
\displaystyle \frac{\lambda_{F_{2}}}{\varphi^{2p_{2}}} + p_{2}
\left[(\zeta-1)\frac{(\varphi^{\prime})^{2}}{\varphi^{2}} +
\frac{\varphi^{\prime \prime}}{\varphi} \right]=\lambda.
  \end{cases}
  \leqno(E-K-II)
\end{equation*}
The last equation implies in particular that $\lambda_{F_{2}}$ is
constant.

Let the system
\begin{equation*}
  \begin{cases}
     (\varphi^{\sigma})^{\prime \prime}=\nu \varphi^{\sigma}\\
     0 < \varphi .
  \end{cases}
  \leqno(\varphi^{\sigma};\nu)
\end{equation*}
where $\nu$ and $\sigma$ are real parameters. All its solutions
$\varphi^\sigma$ have the form
\begin{equation*}
\varphi^\sigma(t)= \begin{cases} \mathcal{A} e^{i\sqrt{-\nu}\, t}
+ \mathcal{B} e^{-i\sqrt{-\nu}\,
t} & \text{ if }\nu < 0, \\
\mathcal{A}t + \mathcal{B} & \text{ if }\nu =0,\\
\mathcal{A} e^{\sqrt{\nu}\, t} + \mathcal{B} e^{-\sqrt{\nu}\, t} &
\text{ if }\nu > 0,
\end{cases}
\end{equation*}
with constants $\mathcal{A}$ and $\mathcal{B}$ such that
$\varphi>0$.

\noindent Furthermore, let the $(\varphi^{\sigma};\nu)$ modified
system
\begin{equation*}
\begin{cases}
(\varphi^{\sigma})^{\prime \prime}= \nu \varphi^{\sigma} \\
[(\varphi^{\sigma})^{\prime}]^{2}=\nu (\varphi^{\sigma})^{2}\\
\varphi >0.
\end{cases}
\leqno(\varphi^{\sigma};\nu;*)
\end{equation*}
Note that $\nu$ must be $>0$. It is easy to verify that all its
solutions are given by
\begin{equation*}\label{}
  \varphi^\sigma(t)= \mathcal{A} e^{\pm \sqrt{\nu}\, t},
\end{equation*}
where $\mathcal{A}$ is a positive constant.

Consider now two cases, namely

\begin{description}
\item[$\underline{\zeta = 0}$] First of all, note that
$\displaystyle p_{2} = - \frac{1}{2} p_{1}$ and $\eta =
\displaystyle \frac{3}{2}p_{1}^{2}$. $\hfill$

\begin{description}
\item[$\underline{\eta = 0}$] Thus, $p_{i}=0$, for all $i$ and $0
= \lambda = \lambda_{F_{2}}$.
Thus the corresponding metric is
\begin{equation*}
-{\rm d}t^{2}+g_{F_{1}}+g_{F_{2}}.
\end{equation*}
\item[$\underline{\eta \neq 0}$] Then $p_{1} \neq 0$, $p_{2} \neq
0$ and $\hfill $
\begin{equation*}
\begin{cases}
\displaystyle \eta \frac{(\varphi^{\prime})^{2}}{\varphi^{2}}
=\lambda\\
\displaystyle
p_{1}\left[-\frac{(\varphi^{\prime})^{2}}{\varphi^{2}} +
\frac{\varphi^{\prime \prime}}{\varphi} \right]=\lambda\\
\displaystyle \frac{\lambda_{F_{2}}}{\varphi^{-p_{1}}} -
\frac{1}{2} p_{1}
\left[-\frac{(\varphi^{\prime})^{2}}{\varphi^{2}} +
\frac{\varphi^{\prime \prime}}{\varphi} \right]=\lambda.
\end{cases}
\leqno(E-K-II i)
\end{equation*}
If
\begin{description}
\item[$\underline{\lambda_{F_{2}}=0}$] then $\lambda = 0$ and
$\varphi$ is constant $\varphi_{0}$. Thus the corresponding
metric is
\begin{equation*}
-{\rm
d}t^{2}+\varphi_{0}^{2p_{1}}g_{F_{1}}+\varphi_{0}^{2p_{2}}g_{F_{2}}.
\end{equation*}
\item[$\underline{\lambda_{F_{2}}\neq0}$] then
$\displaystyle\frac{\lambda_{F_{2}}}{\varphi^{-p_{1}}} =
\frac{3}{2} \lambda$, as consequence $\varphi$ is constant and
considering the system this gives a contradiction. $\hfill $
\end{description}
\end{description}

\item[$\underline{\zeta \neq 0}$] Hence $ \eta \neq 0$ and by
{\it Remark \ref{rem3}} the system reduces to
\begin{equation*}
\begin{cases}
\displaystyle \frac{\zeta^{2}}{\eta} \frac{(
\varphi^{\zeta\frac{\eta}{\zeta^{2}}})^{\prime \prime}}
{\varphi^{\zeta\frac{\eta}{\zeta^{2}}}} = \lambda\\
\displaystyle \frac{p_{1}}{\zeta}
\frac{(\varphi^{\zeta})^{\prime \prime}}{\varphi^{\zeta}}=\lambda \\
\displaystyle \frac{\lambda_{F_{2}}}{\varphi^{2p_{2}}} +
\frac{p_{2}}{\zeta} \frac{(\varphi^{\zeta})^{\prime
\prime}}{\varphi^{\zeta}} =\lambda,
\end{cases}
\leqno(E-K-II ii)
\end{equation*}

\begin{description}
\item[$\underline{\eta=\zeta^{2}}$] So $p_{1}\neq 0$ and either
$p_{2} = 0$ or $p_{2} = -2p_{1}$.
If $\hfill$
\begin{description}
\item[$\underline{\lambda = 0}$] then $ \lambda_{F_{2}}=0$. Thus,
the corresponding metric is $\hfill $
\begin{equation*}
-{\rm d}t^{2}+\varphi^{2p_{1}}
g_{F_{1}}+\varphi^{2p_{2}}g_{F_{2}},
\end{equation*}
where $\varphi$ satisfies $(\varphi^{\zeta};0)$.
\item[$\underline{\lambda \neq 0}$]  then $\zeta=p_{1}$ and
$p_{2}=0$. Hence, by the third equation $\hfill $
$\lambda_{F_{2}}=\lambda$. Thus, the corresponding metric is
\begin{equation*}
-{\rm d}t^{2}+\varphi^{2\zeta} g_{F_{1}}+g_{F_{2}},
\end{equation*}
where $\varphi$ satisfies $(\varphi^{\zeta};\lambda)$.
\end{description}
\item[$\underline{\eta\neq\zeta^{2}}$] Then $p_{2}\neq 0$
and $\displaystyle
\frac{p_{1}}{p_{2}}\frac{\lambda_{F_{2}}}{\varphi^{2p_{2}}}
=\left(\frac{p_{1}}{p_{2}}-1\right)\lambda$.

So, if
\begin{description}
\item[$\underline{\lambda=0}$] then the first equation implies,
$\varphi^{\zeta}=(\mathcal{A}t+\mathcal{B})^{\frac{\zeta^{2}}{\eta}}$
and $\displaystyle(\varphi^{\zeta}) ^{\prime
\prime}=\frac{\zeta^{2}}{\eta}\left(\frac{\zeta^{2}}{\eta}-1\right)
(\mathcal{A}t+\mathcal{B})^{\frac{\zeta^{2}}{\eta}-2}\mathcal{A}^{2}$.

\begin{description}
\item[$\underline{\lambda_{F_{2}}=0}$]then applying the third
equation results $\mathcal{A}=0$, so $\varphi^{\zeta}$ is
constant and $\varphi$ is a positive constant $\varphi_{0}$.
Thus the corresponding metric is
\begin{equation*}
-{\rm
d}t^{2}+\varphi_{0}^{2p_{1}}g_{F_{1}}+\varphi_{0}^{2p_{2}}g_{F_{2}}.
\end{equation*}
\item[$\underline{\lambda_{F_{2}}\neq0}$]then $p_{1}=0$, hence
$\displaystyle p_{2}=\frac{\zeta}{2}$, $\displaystyle
\frac{\eta}{\zeta^{2}}=\frac{1}{2}$. So, by the third equation
$\lambda_{F_{2}}=-\mathcal{A}^{2}<0$. Thus the corresponding
metric is
\begin{equation*}
-{\rm d}t^{2}+g_{F_{1}}+\varphi^{2p_{2}}g_{F_{2}},
\end{equation*}
with $\varphi$ as above.
\end{description}
\item[$\underline{\lambda\neq0}$] then $p_{1}\neq 0$, hence
$\displaystyle \frac{\lambda_{F_{2}}}{\varphi^{2p_{2}}}=\left(1
- \frac{p_{1}}{p_{2}}\right)\lambda$.
\begin{description}
\item[$\underline{\lambda_{F_{2}}=0}$] then $p_{1}=p_{2}$ and the
system can be reduced to
\begin{equation*}
3\frac{(\varphi^{\zeta \frac{1}{3}})^{\prime
\prime}}{\varphi^{\zeta \frac{1}{3}}} =
\frac{1}{3}\frac{(\varphi^{\zeta})^{\prime
\prime}}{\varphi^{\zeta}} = \lambda
\end{equation*}
which is equivalent to the solvable system
$(\varphi^{\zeta};3\lambda;*)$.
%\begin{equation*}
%\begin{cases}
%(\varphi^{\zeta})^{\prime \prime}= 3\lambda \varphi^{\zeta} \\
%[(\varphi^{\zeta})^{\prime}]^{2}=3\lambda (\varphi^{\zeta})^{2}\\
%\varphi >0.
%\end{cases}
%\leqno(*)
%\end{equation*}
Note that $\lambda$ must be $>0$.

\item[$\underline{\lambda_{F_{2}}\neq0}$]then $\varphi$ is
constant and this gives a contradiction. $\hfill $
\end{description}
\end{description}
\end{description}
\end{description}

The table that follows specifies the only possible Einstein
generalized Kasner space-times of Type (II) with the corresponding
parameters. The last column indicates the function $\varphi$ or
the system which it satisfies.

\begin{center}
\small {\begin{tabular}{|c|c|c|c|c|c|c|c|c|c|} \hline $\zeta$ &
$\eta$ & $\frac{\eta}{\zeta^{2}}$ & $\lambda$ &
$\lambda_{F_{2}}$ & $p_{1}$ & $p_{2}$ & metric & $\varphi$  \\
\hline 0 & 0 & - & 0 & 0 & 0 & 0 & $-{\rm d}t^{2}+g_{F_{1}}+g_{F_{2}}$ & -  \\
0 & $\frac{3}{2} p_{1}^{2}\neq 0 $ & - & 0 & 0 & $\neq 0$ &
$-\frac{1}{2} p_{1} $ & $-{\rm
d}t^{2}+\varphi_{0}^{2p_{1}}g_{F_{1}}+\varphi_{0}^{-p_{1}}g_{F_{2}}$
& $\varphi_{0}=cte>0$  \\
0 & $\frac{3}{2} p_{1}^{2} \neq 0 $ & - & $-$ & $\neq 0$ & $\neq
0 $ & $-\frac{1}{2}p_{1} $ &
no metric &   - \\
$\neq 0$& $\zeta^{2} $ & 1 & $ 0 $ & $ 0 $ & $\neq 0$ & $0,
-2p_{1} $ & $-{\rm
d}t^{2}+\varphi^{2p_{1}}g_{F_{1}}+\varphi^{2p_{2}}g_{F_{2}}$ &
$(\varphi^{\zeta};0)$  \\
$\neq 0$& $\zeta^{2} $ & 1 & $\neq 0 $ & $\lambda $ & $\neq 0$ &
0 & $-{\rm d}t^{2}+\varphi^{2p_{1}}g_{F_{1}}+g_{F_{2}}$ &
$(\varphi^{\zeta};\lambda)$   \\
$\neq 0$& $\neq 0$ & $\neq 1$ & 0 & 0 &  $p_{1}$ & $\neq 0 $ &
$-{\rm
d}t^{2}+\varphi_{0}^{2p_{1}}g_{F_{1}}+\varphi_{0}^{2p_{2}}g_{F_{2}}$
& $\varphi_{0}=cte>0$ \\
$\neq 0$& $\neq 0$ & $\neq 1$ & 0 & $< 0$ &  0 & $\neq 0 $ &
$-{\rm d}t^{2}+g_{F_{1}}+\varphi^{2p_{2}}g_{F_{2}}$ &
$(\varphi^{\frac{\eta}{\zeta}};0)$  \\
$\neq 0$& $\neq 0$ & $\neq 1$ & $>0 $ &  0 & $p_{2} $ & $\neq 0
$ &
$-{\rm d}t^{2}+\varphi^{2p_{1}}g_{F_{1}}+\varphi^{2p_{1}}g_{F_{2}}$
& $(\varphi^{\zeta};3\lambda;*)$\\
$\neq 0$& $\neq 0$ & $\neq 1$ & $\neq 0$ &  $\neq 0$ & $p_{1}$ &
$\neq 0 $ &  no metric & -\\
\hline
\end{tabular}}
\end{center}

\begin{center} {\sc Table 1} \end{center}

Note that {\it Corollary \ref{cor3}} cannot be applied in the
situations above.

$\bullet$ {\bf Classification of the Type (II) generalized
Kasner space-times with constant scalar curvature}

Let $M = I \times_{\varphi^{p_{1}}}F_{1}\times_{\varphi^{p_{2}}}
F_{2}$ be a Type (II) generalized Kasner space-time with constant
scalar curvature. Then the parameters introduced before {\it
Proposition \ref{ks-2}} satisfy $\zeta = p_{1}+2p_{2}$, $\eta =
p_{1}^{2}+2p_{2}^{2}$ and
\begin{equation*}
  \tau = 2 \zeta \frac{\varphi^{\prime \prime}}{\varphi} + [(\zeta - 2)\zeta
+\eta]\frac{(\varphi^{\prime})^{2}}{\varphi^{2}}
  + \frac{\tau_{F_{2}}}{\varphi^{2p_{2}}}.
  \leqno(csc-K-II.a)
\end{equation*}
Note that $\tau_{F_{2}}$ must be constant if there exist a
positive solution of $(csc-K-II.a)$ (see also {\it Proposition
\ref{gscamu2}}). We consider two principal cases with different
subcases.

\begin{description}
\item[$\underline{\zeta=0}$] If
\begin{description}
\item[$\underline{\eta =0}$] then $p_{1}=p_{2}=0$, $\tau =
\tau_{F_{2}}$ and the corresponding metric is
\begin{equation*}
-{\rm d}t^{2}+g_{F_{1}}+g_{F_{2}}.
\end{equation*}
\item[$\underline{\eta \neq 0}$] then $\displaystyle
p_{2}=-\frac{1}{2}p_{1}$ and $\displaystyle \eta =
\frac{3}{2}p_{1}^{2}= 6 p_{2}^{2}$. The equation $(csc-K-II.a)$
reduces to \begin{equation*} \tau = \eta
\frac{(\varphi^{\prime})^{2}}{\varphi^{2}} +
\frac{\tau_{F_{2}}}{\varphi^{2p_{2}}}. \leqno(csc-K-II.b)
\end{equation*}
\end{description}
\item[$\underline{\zeta\neq 0}$] implies $ \eta\neq 0$ and
considering
$0<u=(\varphi^{\zeta})^{\frac{1+\frac{\eta}{\zeta^{2}}}{2}}$,
{\it Corollary \ref{cor13}} arises the relation
\begin{equation*}\label{}
\displaystyle -\frac{4}{1+\displaystyle \frac{\eta}{\zeta^{2}}}
u^{\prime \prime}= -\tau u+
\tau_{F_{2}}u^{1-\frac{4}{1+\frac{\eta}
{\zeta^{2}}}\frac{p_{2}}{\zeta}}. \leqno(csc-K-II.c)
\end{equation*}
\begin{description}
\item[$\underline{\eta=\zeta^{2}}$]Then $p_{1}\neq 0$, either
$p_{2} = 0$ or $p_{2}=-2p_{1}$, and $u = \varphi^{\zeta}$.
\begin{description}
\item[$\underline{\tau_{F_{2}}  =0}$] So the equation reduces to
\begin{equation*}
-2u^{\prime \prime}=-\tau u.
\end{equation*}
\item[$\underline{\tau_{F_{2}}  \neq 0}$] If
\begin{description}
\item[$\underline{p_{2} = 0}$] the equation reduces to
\begin{equation*}
-2u^{\prime \prime}=-(\tau  - \tau_{F_{2}})u.
\end{equation*}
\item[$\underline{p_{2}=-2p_{1}}$]
\begin{equation*}
-2u^{\prime \prime}=-\tau u + \tau_{F_{2}}u^{-\frac{1}{3}}.
%\leqno(3)
\end{equation*}
\end{description}
\end{description}
\item[$\underline{\eta\neq\zeta^{2}}$] Then $p_{2}\neq 0$ and $\displaystyle
\frac{\eta}{\zeta^{2}}\ge\frac{1}{3}$.
\begin{description}
\item[$\underline{\tau_{F_{2}}=0}$]
\begin{equation} \label{IIs}
\displaystyle -\frac{4}{1+\displaystyle \frac{\eta}{\zeta^{2}}}
u^{\prime \prime}= -\tau u
\end{equation}
\item[$\underline{\tau_{F_{2}}\neq 0}$]
\begin{equation*}\label{}
\displaystyle -\frac{4}{1+\displaystyle \frac{\eta}{\zeta^{2}}}
u^{\prime \prime}= -\tau u+
\tau_{F_{2}}u^{1-\frac{4}{1+\frac{\eta}{\zeta^{2}}}\frac{p_{2}}{\zeta}}.
\leqno(csc-K-II.c)
\end{equation*}

Note that a particular subcase is $\displaystyle
\frac{\eta}{\zeta^{2}}=\frac{1}{3}.$ In fact, in this case,
$p_{1}=p_{2}=\displaystyle \frac{\zeta}{3}$ (see {\it Remark
\ref{rem-i}}) and the latter equation reduces to the
non-homogeneous linear ordinary differential equation
\begin{equation*} \label{}
\displaystyle -3u^{\prime \prime}= -\tau u+ \tau_{F_{2}}.
\end{equation*}
\end{description}
\end{description}
\end{description}

Synthetically, remembering that in each case the corresponding
metric may be written as $ -{\rm d}t^{2}+\varphi^{2p_{1}}
g_{F_{1}}+\varphi^{2p_{2}} g_{F_{2}}$, we find that the only
possibilities to have constant scalar curvature in a generalized
Kasner space-time of type (II) are generated by

\begin{center}{
\begin{tabular}{|c|c|c|c|c|c|c|c|c|c|}
\hline $\zeta$ & $\eta$ & $\frac{\eta}{\zeta^{2}}$   &
$\tau_{F_{2}}$ & $p_{1}$ & $p_{2}$  & $\varphi$ eq.  \\ \hline
0 & 0 & -   & $\tau_{F_{2}}$ & 0 & 0 &    $\tau = \tau_{F_{2}}$ \\
0 & $\frac{3}{2}p_{1}^{2}$ & -    & 0 & $\neq 0$ &
$-\frac{1}{2}p_{1}$ & $\tau=\eta\frac{(\varphi^{\prime})^{2}}
{\varphi ^{2}}$ \\
0 & $\frac{3}{2}p_{1}^{2}$ & - & $ \neq 0$ & $\neq 0$ &
$-\frac{1}{2}p_{1}$ & $\tau=\eta\frac{(\varphi ^{\prime})
^{2}}{\varphi ^{2}}+ \frac{\tau_{F_{2}}}{\varphi ^{2p_{2}}}$ \\
$\zeta \neq 0$ & $\zeta^{2} $ & 1    & 0 & $\neq 0$ & 0 &
$-2u^{\prime \prime}=-\tau u  ;  \hfill u=\varphi^{\zeta}$\\
$\zeta \neq 0$ & $\zeta^{2} $ & 1    & 0 & $\neq 0$ & $-2p_{1}$
& $-2u^{\prime \prime}=-\tau u  ; \hfill u=\varphi^{\zeta}$\\
$\zeta \neq 0$ & $\zeta^{2} $ & 1    & $\neq 0$ & $\neq 0$ & $0$
& $-2u^{\prime \prime}=-(\tau -\tau_{F_{2}})  ; \hfill u=\varphi^{\zeta}$\\
$\zeta \neq 0$ & $\zeta^{2} $ & 1    & $\neq 0$ & $\neq 0$ &
$-2p_{1}$ & $-2u^{\prime \prime}=-\tau u +
\tau_{F_{2}}u^{-\frac{1}{3}} ; \hfill
u=\varphi^{\zeta}$ \\
$\zeta \neq 0$ & $\eta \neq 0$ & $\neq 1$    & 0 & $p_{1}$ & $\neq
0$ & (\ref{IIs}); \hfill
$u=(\varphi^{\zeta})^{\frac{1+\frac{\eta}{\zeta^{2}}}{2}}$ \\
$\zeta \neq 0$ & $\eta \neq 0$ & $\neq 1,\frac{1}{3}$    & $\neq
0$ & $p_{1}$ & $\neq 0$ & $(csc-K-II.c); \hfill
u=(\varphi^{\zeta})^{\frac{1+\frac{\eta}{\zeta^{2}}}{2}}$ \\
$\zeta \neq 0$ & $\frac{\zeta^{2}}{3} $ & $\frac{1}{3}$ & $\neq
0$ & $\frac{\zeta}{3}$ & $\frac{\zeta}{3}$ & $-3u^{\prime
\prime} = -\tau u +\tau_{F_{2}}; \hfill u =
\varphi^{\frac{2}{3}\zeta}$ \\
\hline
\end{tabular}}
\end{center}

\begin{center} {\sc Table 2} \end{center}

\noindent where the conditions for $\tau$ must be imposed by the
existence of positive solutions of the ordinary differential
equations of the last column, on the corresponding interval $I$.

\subsection{Type (III)} \label{chp5c}

Let $M=I \times {}_{b_1}F_1 \times {}_{b_2}F_2 \times {}_{b_3}F_3$
be a type (III) warped product space-time with metric $g=-{\rm
d}t^2 \oplus b_{1}^{2}g_{F_1} \oplus b_{2}^{2}g_{F_2} \oplus
b_{3}^2g_{F_3}.$ Then the scalar curvature $\tau$ of $(M,g)$ is
given as

$$\tau=2\left( \frac{b_1^{\prime \prime}}{b_1}+ \frac{b_2^{\prime
\prime}}{b_2}+\frac{b_3^{\prime \prime}}{b_3}+ \frac{b_1^\prime
b_2^\prime}{b_1 b_2}+\frac{b_2^\prime b_3^\prime}{b_2
b_3}+\frac{b_1^\prime b_3^\prime}{b_1 b_3} \right).$$

Note that $\tau_{F_i}=0,$ since ${\rm dim}(F_i)=1,$ for any $i
\in \{1,2,3\}.$

If $v_i$ is a vector field on $F_i,$ for any $i \in \{1,2,3\}$ and
$\displaystyle{\overline x=\frac{\partial}{\partial
t}+v_1+v_2+v_3},$ then \begin{eqnarray*} \Ric(\overline
x,\overline x)& = & \left( b_1 b_1^{\prime \prime} + \frac{b_1
b_1^\prime b_2^\prime}{b_2} + \frac{b_1 b_1^\prime
b_3^\prime}{b_3} \right) g_{F_1}(v_1,v_1) \\
& + & \left( b_2 b_2^{\prime \prime}+ \frac{b_2 b_2^\prime
b_1^\prime}{b_1} + \frac{b_2 b_2^\prime
b_3^\prime}{b_3} \right) g_{F_2}(v_2,v_2) \\
& + & \left( b_3 b_3^{\prime \prime}+ \frac{b_3 b_3^\prime
b_1^\prime}{b_1}+\frac{b_3 b_3^\prime
b_2^\prime}{b_2} \right) g_{F_3}(v_3,v_3) \\
& - & \frac{b_1^{\prime \prime}}{b_1}- \frac{b_2^{\prime
\prime}}{b_2}- \frac{b_3^{\prime \prime}}{b_3}
\end{eqnarray*}

Note that $\Ric_{F_i} \equiv 0,$ since ${\rm dim}(F_i)=1,$ for
any $i \in \{1,2,3\}.$

$\bullet$ {\bf Classification of Einstein Type (III) generalized
Kasner space-times}

Let $M = I \times_{\varphi^{p_{1}}}F_{1}\times_{\varphi^{p_{2}}}
F_{2}\times_{\varphi^{p_{3}}}F_{3}$ be an Einstein Type (III)
generalized Kasner space-time. Then the parameters introduced
before {\it Proposition \ref{ks-1}} satisfy $\zeta =
p_{1}+p_{2}+p_{3}$, $\eta = p_{1}^{2}+p_{2}^{2}+p_{3}^{2}$. Hence
the latter arises
\begin{equation*}
\begin{cases}
\displaystyle (\eta -
\zeta)\frac{(\varphi^{\prime})^{2}}{\varphi^{2}} + \zeta
\frac{\varphi^{\prime \prime}}{\varphi} =
\lambda \hfill  \\
\displaystyle
p_{1}\left[(\zeta-1)\frac{(\varphi^{\prime})^{2}}{\varphi^{2}} +
\frac{\varphi^{\prime \prime}}{\varphi} \right]=\lambda\\
\displaystyle
p_{2}\left[(\zeta-1)\frac{(\varphi^{\prime})^{2}}{\varphi^{2}} +
\frac{\varphi^{\prime \prime}}{\varphi}
\right]=\lambda \\
\displaystyle
p_{3}\left[(\zeta-1)\frac{(\varphi^{\prime})^{2}}{\varphi^{2}} +
\frac{\varphi^{\prime \prime}}{\varphi} \right]=\lambda.
\end{cases} \leqno(E-K-III)
\end{equation*}
Note that adding the last three equations, there results
\begin{equation} \label{IIIe}
\displaystyle
\zeta\left[(\zeta-1)\frac{(\varphi^{\prime})^{2}}{\varphi^{2}} +
\frac{\varphi^{\prime \prime}}{\varphi} \right]=3 \lambda
\end{equation}

%Let the system
%\begin{equation*}
%\begin{cases}
%(\varphi^{\sigma})^{\prime \prime}=\nu \varphi^{\sigma}\\
%0 < \varphi.
%\end{cases}
%\leqno(\varphi^{\sigma};\nu),
%\end{equation*}
%where $\nu$ and $\sigma$ are real parameters.

Consider now two cases, namely
\begin{description}
\item[$\underline{\zeta = 0}$] Then applying (\ref{IIIe}), we
obtain $\lambda=0.$
\begin{description}
\item[$\underline{\eta = 0}$] Thus $p_{i}=0$ for all $i$.
Hence the corresponding metric is
\begin{equation*}
-{\rm d}t^{2}+g_{F_{1}}+g_{F_{2}}+g_{F_{3}}.
\end{equation*}

\item[$\underline{\eta \neq 0}$] the system reduces to
\begin{equation*}
\begin{cases}
\displaystyle \eta \frac{(\varphi^{\prime})^{2}}{\varphi^{2}}=0\\
\displaystyle
p_{i}\left[-\frac{(\varphi^{\prime})^{2}}{\varphi^{2}} +
\frac{\varphi^{\prime \prime}}{\varphi} \right]=0 \textrm{ for
all } i=1,2,3.
\end{cases}
\leqno(E-K-III i)
\end{equation*}
then $\varphi$ is constant $\varphi_{0}$. Thus the corresponding
metric is
\begin{equation*}
-{\rm d}t^{2}+\varphi_{0}^{2p_{1}}g_{F_{1}}+\varphi_{0}^{2p_{2}}
g_{F_{2}}+\varphi_{0}^{2p_{3}}g_{F_{3}}
\end{equation*}
\end{description}

\item[$\underline{\zeta \neq 0}$] Thus $ \eta \neq 0$ and by
{\it Remark \ref{rem3}} the system reduces to
\begin{equation*}
\begin{cases}
\displaystyle \frac{\zeta^{2}}{\eta} \frac{(
\varphi^{\zeta\frac{\eta}{\zeta^{2}}})^{\prime \prime}}
{\varphi^{\zeta\frac{\eta}{\zeta^{2}}}} = \lambda\\
\displaystyle \frac{p_{i}}{\zeta}
\frac{(\varphi^{\zeta})^{\prime \prime}}{\varphi^{\zeta}}
=\lambda \textrm{ for all } i=1,2,3.
\end{cases} \leqno(E-K-III ii)
\end{equation*}
Adding the last three equations in $(E-K-III ii)$, we obtain that
\begin{equation} \label{IIIe2}
\displaystyle \frac{(\varphi^{\zeta})^{\prime
\prime}}{\varphi^{\zeta}} =3 \lambda
\end{equation}

\begin{description}
\item[$\underline{\eta=\zeta^{2}}$] Then (\ref{IIIe}) and
(\ref{IIIe2}), give $\lambda = 0$. Thus, the corresponding metric
is
\begin{equation*}
-{\rm d}t^{2}+\varphi^{2p_{1}}
g_{F_{1}}+\varphi^{2p_{2}}g_{F_{2}} +\varphi^{2p_{3}}g_{F_{3}},
\end{equation*}
where $\varphi$ satisfies $(\varphi^{\zeta};0)$.

\item[$\underline{\eta\neq\zeta^{2}}$] Then at least two $p_{i}$'s
are $\neq 0$. So if
\begin{description}
\item[$\underline{\lambda=0}$] then the first equation implies
$\varphi^{\zeta}=(\mathcal{A}t+\mathcal{B})^{\frac{\zeta^{2}}{\eta}}$
and $\displaystyle(\varphi^{\zeta}) ^{\prime
\prime}=\frac{\zeta^{2}}{\eta}\left(\frac{\zeta^{2}}{\eta}-1\right)
(\mathcal{A}t+\mathcal{B})^{\frac{\zeta^{2}}{\eta}-2}\mathcal{A}^{2}$.
Then by (\ref{IIIe2}) results $\mathcal{A}=0$, so
$\varphi^{\zeta}$ is constant and $\varphi$ is a positive constant
$\varphi_{0}$. Thus the corresponding metric is
\begin{equation*}
-{\rm d}t^{2}+\varphi_{0}^{2p_{1}}g_{F_{1}}+\varphi_{0}^{2p_{2}}
g_{F_{2}}+\varphi_{0}^{2p_{3}}g_{F_{3}}.
\end{equation*}

\item[$\underline{\lambda\neq0}$] then all $p_{i}$'s are $\neq 0$
and all of them are equals, so that $\displaystyle
p_{1}=p_{2}=p_{3}=\frac{\zeta}{3}$. So $\displaystyle
\eta=\frac{\zeta^{2}}{3}$ and $\displaystyle
\frac{\eta}{\zeta^{2}}=\frac{1}{3}$. Thus the system reduces to
\begin{equation*}
3\frac{(\varphi^{\zeta \frac{1}{3}})^{\prime
\prime}}{\varphi^{\zeta \frac{1}{3}}} =
\frac{1}{3}\frac{(\varphi^{\zeta})^{\prime
\prime}}{\varphi^{\zeta}} = \lambda,
\end{equation*}
which is equivalent to the solvable system
$(\varphi^{\zeta};3\lambda;*)$.
%\begin{equation*}
%\begin{cases}
%(\varphi^{\zeta})^{\prime \prime}= 3\lambda \varphi^{\zeta} \\
%[(\varphi^{\zeta})^{\prime}]^{2}=3\lambda (\varphi^{\zeta})^{2}\\
%\varphi >0.
%\end{cases}
%\leqno(**)
%\end{equation*}
Note that $\lambda$ must be $>0$.
\end{description}
\end{description}

The table that follows specifies the only possible Einstein
generalized Kasner space-times of type (III)  with the
corresponding parameters. Like for the table of Type (II), the
last column indicates the function $\varphi$ or the system which
it satisfies.

This example may be easily generalized to the situation all
the $F_{i}$'s are Ricci flat, considering $S=\sum_{i=1}^{m}
s_{i}
>1$ instead of $3$.
\end{description}

\begin{center}
\small
\begin{tabular}{|c|c|c|c|c|c|c|c|c|c|}
\hline $\zeta$ & $\eta$ & $\frac{\eta}{\zeta^{2}}$ & $\lambda$ &
$p_{1}$ & $p_{2}$ & $p_{3}$ & metric & $\varphi$  \\ \hline 0 &
0 & - & 0 & 0 & 0 & 0 &
$-{\rm d}t^{2}+g_{F_{1}}+g_{F_{2}}+g_{F_{3}}$ & - \\
0 & $\neq 0$ & - & 0 & $p_{1}$ & $p_{2}$ & $p_{3}$ & $-{\rm
d}t^{2}+\varphi_{0}^{2p_{1}}g_{F_{1}}+\varphi_{0}^{2p_{2}}
g_{F_{2}}+\varphi_{0}^{2p_{3}}g_{F_{3}}$
& $\varphi_{0}=cte>0$  \\
$\neq 0$& $\zeta^{2}$ & 1 & $0$ & $p_{1}$ & $p_{2}$ & $p_{3}$ &
$-{\rm d}t^{2}+\varphi^{2p_{1}}g_{F_{1}}+
\varphi^{2p_{2}}g_{F_{2}}+ \varphi^{2p_{3}}g_{F_{3}}$
& $(\varphi^{\zeta};0)$ \\
$\neq 0$& $\neq 0$ & $\neq 1$ & 0 & $p_{1}$ & $p_{2}$ & $p_{3}$
& $-{\rm d}t^{2}+\varphi_{0}^{2p_{1}}g_{F_{1}}+
\varphi_{0}^{2p_{2}}g_{F_{2}}+ \varphi_{0}^{2p_{3}}g_{F_{3}}$
& $\varphi_{0}=cte>0$ \\
$\neq 0$& $\neq 0$ & $\neq 1$ & $> 0$ & $p_{1}$ & $p_{1}$ &
$p_{1}$ & $-{\rm d}t^{2}+\varphi^{2p_{1}}g_{F_{1}}+
\varphi^{2p_{1}}g_{F_{2}}+ \varphi^{2p_{1}}g_{F_{3}}$
& $(\varphi^{\zeta};3\lambda;*)$ \\
\hline
\end{tabular}
\end{center}

\begin{center} {\sc Table 3} \end{center}

$\bullet$ {\bf Classification of Type (III) generalized Kasner
space-times with constant scalar curvature}

Let $M = I \times_{\varphi^{p_{1}}}F_{1}\times_{\varphi^{p_{2}}}
F_{2}\times_{\varphi^{p_{3}}}F_{3}$ be a Type (III) generalized
Kasner manifold with constant scalar curvature. Then the
parameters introduced before {\it Proposition \ref{ks-2}} satisfy
$\zeta = p_{1}+p_{2}+p_{3}$, $\eta =
p_{1}^{2}+p_{2}^{2}+p_{3}^{2}$. Thus, this case is already
included in the analysis of {\it Example \ref{ex1}}.
%
%Hence the latter arises
%\begin{equation*}
%\tau = 2 \zeta \frac{\varphi^{\prime \prime}}{\varphi} + [(\zeta
%- 2)\zeta +\eta]\frac{(\varphi^{\prime})^{2}}{\varphi^{2}}.
%\leqno(csc-K-III)
%\end{equation*}
%
%\begin{description}
%\item[$\underline{\zeta = 0}$] The corresponding problem is
%\begin{equation*}
%\begin{cases}
%(\varphi^{\prime})^{2}=\displaystyle \frac{\tau}{\eta} \varphi^{2}\\
%\varphi > 0,
%\end{cases}
%\leqno(csc-K-III-a)
%\end{equation*}
%which is solvable. Compare with {\it Remark \ref{rem3}}.
%\item[$\underline{\zeta \neq 0}$] It is a particular case of
%{\it Example \ref{ex1}}.
%\end{description}

We will close this section by an example and the following comment
which gives some preliminary ideas about our future plans on this
topic (see also the last section for details).
%
%\notemarg
%

\begin{ex} \label{exsph} Let $M=I \times_{\varphi^{p_{1}}} \mathbb S_{3}
\times_{\varphi^{p_{2}}} \mathbb S_{2}$ be a generalized Kasner
manifold with constant scalar curvature. Then the parameters
introduced before {\it Proposition \ref{ks-2}} are given by $\zeta
= 3p_{1}+2p_{2}$, $\eta = 3p_{1}^{2}+2p_{2}^{2}$. Consider now
$p_{1}=1$ and $p_{2}=-1$, then $\zeta=1$ and $\eta=5$. Hence,
applying {\it Corollary \ref{cor13}} the latter conditions arise
for $u=\varphi^{3}$ the problem
\begin{equation} \label{IIIs}
\begin{cases}
\displaystyle -\frac{2}{3} u^{\prime \prime} + \tau u =
\tau_{\mathbb S_{3}} u^{1-\frac{2}{3}} + \tau_{\mathbb S_{2}}
u^{1+\frac{2}{3}}, \\
u>0,
\end{cases}
\end{equation}
where $\tau_{\mathbb S_{3}}, \tau_{\mathbb S_{2}}>0$ are the
constant scalar curvatures of the corresponding spheres. Note that
the equation in (\ref{IIIs}) has always the constant solution zero
and there exists $\tau_{1}>0$ such that for $\tau = \tau_{1}$
there is only one constant solution  of (\ref{IIIs}) and for any
$\tau > \tau_{1}$ there are two constant solutions of
(\ref{IIIs}), so that there exists a range of $\tau$'s,
$(\tau_{1},+\infty)$, where the problem (\ref{IIIs}) has
multiplicity of solutions; while there is no constant solutions
when $\tau < \tau_{1}$.

On the other hand, as in {\it Example \ref{exsph}}, considering
$\mathbb S_{3}$ instead of $\mathbb S_{2}$ with the same values of
$p_{1}$ and $p_{2}$, i.e., $M = I \times_{\varphi^{p_{1}}} \mathbb
S_{3}\times_{\varphi^{p_{2}}} \mathbb S_{3}$, results $\zeta =
3p_{1}+3p_{2}=0$, $\eta = 3p_{1}^{2}+3p_{2}^{2}=6$. Hence,
applying {\it Proposition \ref{ks-2}} the latter conditions arise
the problem
\begin{equation} \label{IIIs2}
\begin{cases}
-6(\varphi^{\prime})^{2}= - \tau \varphi^{2} + \tau_{S_{3}} +
\tau_{\mathbb S_{3}} \varphi^{4}\\
\varphi > 0 .
\end{cases}
\end{equation}
The equation in (\ref{IIIs2}) does not have the constant solution
zero. Furthermore there is no constant solution of (\ref{IIIs2})
if $\tau < 2 \tau_{\mathbb S_{3}} $, there is only one constant
solution of (\ref{IIIs2}) if $\tau = 2 \tau_{\mathbb S_{3}} $ and
two constant solutions of (\ref{IIIs2}) if $\tau > 2 \tau_{\mathbb
S_{3}}.$
\end{ex}

The cases considered above are just some examples for the
different type of differential equations involved in the problem
of constant scalar curvature when the dimensions, curvatures and
parameters have different values. In a future article, we deal
with the problem of constant scalar curvature of a
pseudo-Riemannian generalized Kasner manifolds with a base of
dimension greater than or equal to $1$. This problem carries to
nonlinear partial differential equations with concave-convex
nonlinearities like in (\ref{IIIs}), among others. Nonlinear
elliptic problems with such nonlinearities have been extensively
studied in bounded domains of $\mathbb{R}^{n}$, after the central
article of Ambrosetti, Brezis and Cerami \cite{ABC}, in which the
authors studied the problem of multiplicity of solutions under
Dirichlet conditions. The problem of constant scalar curvature in
a generalized Kasner manifolds with base of dimension greater than
or equal to $1$ is one of the first examples where those
nonlinearities appear naturally. Another related case is the base
conformal warped products, studied in \cite{SBC}.

\section{BTZ (2+1) Black Hole Solutions} \label{chp6}

Now we consider BTZ (2+1)-Black Hole Solutions and give another
characterization of (BTZ) black hole solutions mentioned in {\it
Section 2} (for further details see \cite{BHTZ,BTZ,HCP,MTZ}) in
order to apply the results obtained in this paper.

All the cases considered in \cite{HCP}, can be obtained applying
the {\it formal} approach that follows. By considering the
corresponding square lapse function $N^2,$ the related
3-dimensional, $(2+1)$-space-time model can be expressed as a
$(2+1)$ multiply generalized Robertson-Walker space-time, i.e.,
\begin{equation}\label{eq: 2+1 btz-ds metric}
{\rm d}s^{2}=-{\rm d}t^{2} + b_{1}^{2}(t) {\rm d}x^{2} +
b_{2}^{2}(t) {\rm d}\phi^{2},
\end{equation}
where
\begin{equation}\label{eq: 2+1 btz-ds warping functions-0}
  \begin{cases}
    b_{1}(t) =  N(F^{-1}(t))\\
    b_{2}(t) =  F^{-1}(t),
  \end{cases}
\end{equation}
with
\begin{equation}\label{eq: F(r)}
  F(r)=\int_{a}^{r} \frac{1}{N(\mu)} {\rm d}\mu
\end{equation}
and $F^{-1}$ the inverse function of $F$ (assuming that there
exists) and $a$ is an appropriate  constant that is most of the
time related to the event horizon.

Recalling
\begin{equation}\label{eq:derivative of the inverse}
  1 = \frac{{\rm d}}{{\rm d}t} (F \circ F^{-1})(t)=
  F^{\prime}(F^{-1}(t))(F^{-1})^{\prime}(t),
\end{equation}
we obtain the following properties by applying the chain rule.
Here, note that all the functions depend on the variable $t$ and
the derivatives are taken with respect to the corresponding
arguments.
\begin{itemize}
  \item $b_{1} = N(b_{2})$
  \item $b_{2}^{\prime}= N(F^{-1})=b_{1}$
  \item $b_{2}^{\prime}= N(b_{2})$
  \item $b_{1}^{\prime}=b_{2}^{\prime \prime}=
  N^{\prime}(b_{2})b_{2}^{\prime}=N^{\prime}(b_{2})b_{1}$
  \item $b_{1}^{\prime \prime}= N^{\prime \prime}(b_{2})
  b_{2}^{\prime}b_{1}+N^{\prime}(b_{2})b_{1}^{\prime}
  = N^{\prime \prime}(b_{2})b_{1}^{2}+(N^{\prime}(b_{2}))^{2}b_{1}$.
\end{itemize}
 Thus,
\begin{equation}\label{eq: 2+1 btz-ds warping functions}
  \begin{array}{rrl}
    \bullet & \displaystyle \frac{b_{1}^{\prime \prime}}{b_{1}} &=
N^{\prime \prime}(b_{2})b_{1} +(N^{\prime}(b_{2}))^{2} \\
     & &= N^{\prime \prime}(b_{2})N(b_{2}) +(N^{\prime}(b_{2})) ^{2} \\
     & &= (N^{\prime}N) ^{\prime} (b_{2})\\
     & &= \displaystyle \frac{1}{2} (N^{2}) ^{\prime \prime} (b_{2})\\
    \bullet & \displaystyle \frac{b_{2}^{\prime \prime}}{b_{2}} &=
  N^{\prime}(b_{2})\displaystyle \frac{N(b_{2})}{b_{2}} \\
     & & = \displaystyle \frac{1}{2} \frac{(N^{2}) ^{\prime}
(b_{2})}{b_{2}}\\
    \bullet &  \displaystyle
  \frac{b_{1}^{\prime}}{b_{1}}\frac{b_{2}^{\prime}}{b_{2}}&=
  \displaystyle \frac{b_{2}^{\prime \prime}}{b_{2}}
  \end{array}
\end{equation}
\bigskip

On the other hand, by {\it Corollary \ref{gscamu1}} applied to the
metric \eqref{eq: 2+1 btz-ds metric}, with $s_{1}=s_{2}=1.$ The
scalar curvature of the corresponding space-time is given by
\begin{equation}\label{eq:btz-ds scalar curvature}
  \begin{array}{lll}
  \tau &=& 2 \left(\displaystyle \frac{b_{1}^{\prime \prime}}{b_{1}} +
\frac{b_{2}^{\prime \prime}}{b_{2}} +
  \frac{b_{1}^{\prime}}{b_{1}}\frac{b_{2}^{\prime}}
  {b_{2}}\right)\\
       &=& 2 \left(\displaystyle \frac{b_{1}^{\prime \prime}}
       {b_{1}} +  2
\frac{b_{2}^{\prime \prime}}{b_{2}} \right)\\
    &=& 2 \left( N^{\prime \prime}(b_{2})N(b_{2}) +
    (N^{\prime}(b_{2})) ^{2} +
    2 N^{\prime}(b_{2})\displaystyle \frac{N(b_{2})}{b_{2}}
    \right)\\
    &=& (N^{2}) ^{\prime \prime} (b_{2})+ 2 \displaystyle
    \frac{(N^{2}) ^{\prime}
    (b_{2})}{b_{2}}
  \end{array}
\end{equation}
Note that, the latter is an expression of the scalar curvature
as an operator in the square lapse function. Remember that
$b_{2} = F^{-1}$.

About the Ricci tensor, applying our {\it Proposition
\ref{gpmriki}} and {\it Theorem \ref{gric-fe}} and by
considering again $s_1=s_2=1$, {\it Theorem \ref{gric-fe}} says
that the metric \eqref{eq: 2+1 btz-ds metric} is Einstein with
$\lambda $ if and only if
\begin{equation}\label{eq: 2+1 btz-ds Einstein}
  \begin{cases}
    \displaystyle \frac{b_{1}^{\prime \prime}}{b_{1}} +
    \frac{b_{2}^{\prime \prime}}{b_{2}}
    = \lambda\\
    \displaystyle \frac{b_{1}^{\prime \prime}}{b_{1}} +
    \frac{b_{1}^{\prime}}{b_{1}}
    \frac{b_{2}^{\prime}}{b_{2}}= \lambda \\
    \displaystyle \frac{b_{2}^{\prime \prime}}{b_{2}} +
    \frac{b_{2}^{\prime}}{b_{2}}
    \frac{b_{1}^{\prime}}{b_{1}}=
    \lambda.
  \end{cases}
\end{equation}
On the other hand by making use of \eqref{eq: 2+1 btz-ds warping
functions}, the system \eqref{eq: 2+1 btz-ds Einstein} is
equivalent to (all the functions are evaluated in $r=b_{2}$)

\begin{equation}\label{eq: 2+1 btz-ds Einstein-1}
  \begin{cases}
   (N^{2}) ^{\prime \prime} +  \displaystyle
   \frac{(N^{2})^{\prime}}{r}= 2 \lambda\\
   (N^{2}) ^{\prime \prime} +  \displaystyle
   \frac{(N^{2})^{\prime}}{r}= 2 \lambda\\
   \displaystyle \frac{(N^{2}) ^{\prime}}{r} = \lambda,
  \end{cases}
\end{equation}
or moreover to the following
\begin{equation}\label{eq: 2+1 btz-ds Einstein-2}
  \begin{cases}
  (N^{2})^{\prime \prime}+ \displaystyle
  \frac{(N^{2})^{\prime}}{r} = 2\lambda\\
  \displaystyle \frac{(N^{2}) ^{\prime}}{r} = \lambda.
  \end{cases}
\end{equation}
Thus, we have
\begin{equation}\label{eq: 2+1 btz-ds Einstein-3}
  (N^{2})^{\prime \prime} = \lambda.
\end{equation}
Hence,
\begin{equation}\label{eq: 2+1 btz-ds Einstein-4}
  N^{2}(r) = \frac{\lambda}{2} r^{2} + c_{1} r + c_{2},
\end{equation}
with $c_{1}$ and $c_{2}$ suitable constants. But, since $(N^{2})
^{\prime} (r) = \lambda r + c _{1}$, the second equation of
\eqref{eq: 2+1 btz-ds Einstein-2} is verified if and only if
$c_{1}=0$. So, we have proved the following results.

\begin{prp} \label{btz-1}
Suppose that we have a $(2+1)$-Lorentzian multiply warped product
with the metric given by \eqref{eq: 2+1 btz-ds metric}, where
$b_{1}$ and $b_{2}$ satisfying both \eqref{eq: 2+1 btz-ds warping
functions-0} and \eqref{eq: F(r)}. The space-time is Einstein with
Ricci curvature $\lambda$ if and only if the square lapse function
$N^{2}$ satisfies \eqref{eq: 2+1 btz-ds Einstein-4}, with
$c_{1}=0$ and a suitable constant $c_{2}.$
\end{prp}

Notice that the static (BTZ) and the static (dS) black hole
solutions considered in \cite{HCP} satisfy {\it Proposition
\ref{btz-1}}. Thus they are Einstein multiply warped product
space-times.

\begin{rem} \label{btz-2} Remark that if $N^{2}$ satisfies
\eqref{eq: 2+1 btz-ds Einstein-4} with $c_{1}=0$, then an
application of \eqref{eq:btz-ds scalar curvature} gives the
constancy of the scalar curvature $\tau = 3\lambda $, as desired.
Note that this result agrees with the ones obtained in \cite{HCP}.
\end{rem}

Furthermore, the following just follows from the solution of the
involved second order linear ordinary differential equation
arisen by the expression \eqref{eq:btz-ds scalar curvature}.

\begin{prp} \label{btz-3} Suppose that there is a
$(2+1)$-Lorentzian multiply warped product with the metric given
by \eqref{eq: 2+1 btz-ds metric}, where $b_{1}$ and $b_{2}$
verifying \eqref{eq: 2+1 btz-ds warping functions-0} and
\eqref{eq: F(r)}. The space-time has constant scalar curvature
$\tau = \lambda$ if and only if the square lapse function $N^{2}$
has the form
\begin{equation}\label{eq:btz-ds constant scalar curvature}
 N^{2}(r)=-c_{1}\frac{1}{r} + \frac{\lambda}{6} r^{2} + c_{2},
\end{equation}
with suitable constants $c_{1}$ and $c_{2}.$
\end{prp}

Note that {\it Proposition \ref{btz-3}} agrees with {\it Remark
\ref{btz-2}}.

\section{Conclusions} \label{chp7}
Now, we would like to summarize the content of the paper and to
make some concluding remarks. In a brief, we studied expressions
that relate the Ricci (respectively scalar) curvature of a
multiply warped product with the Ricci (respectively scalar)
curvatures of its base and fibers as well as warping functions.

By using expressions obtained in the paper, we proved necessary
and sufficient conditions for a multiply generalized
Robertson-Walker space-time to be Einstein or to have constant
scalar curvature.

Furthermore, we introduced and considered a kind of generalization
of Kasner space-times, which is closely related to recent
applications in cosmology where metrics of the form
\begin{equation}\label{eq:g-k1}
  {\rm d}s^2 = - {\rm d}t^2 + \sum_{i=1}^{k} e^{2 \alpha_i}  {\rm d}x_i^2,
  \textrm{ with } \alpha_i = \alpha_i(t),
\end{equation}
are frequently considered (see
\cite{{Halpern00},{Socorro-Villanueva-Pimentel96}}; for other
recent topics concerned Kasner type metrics see for instance
\cite{{Dabrowski99},{Frolov01},{Ivashchuk-Singleton04},{JMS},{Maceda04},{Papadopoulos04},{van
den Hoogen-Horne04},{Vladimirov-Kokarev02}}).
%[26,39,46,47,60,65,74,\cite{Vladimirov-Kokarev02}]).
If each warping function $e^{2 \alpha_i}$ is expressed as
\begin{equation}\label{eq:g-k2}
 e^{2 \alpha_i}= \varphi_i^{2 p_i}, \textrm{ with }
 \varphi_i=e^{\frac{\alpha_i}{p_i}},
\end{equation}
for suitable $p_i$'s, then \eqref{eq:g-k1} takes the form
\begin{equation}\label{eq:g-k3}
   {\rm d}s^2 = -{\rm d}t^2 + \sum_{i=1}^{k} \varphi_i^{2 p_i} {\rm d}x_i^2.
\end{equation}

\noindent Our generalization of Kasner space-times corresponds
exactly to the case in which the $\varphi_i$'s are independent of
$i$. More explicitly, $\alpha_i = p_i \alpha$ in \textit{Equation}
\eqref{eq:g-k2}, with $\alpha=\alpha(t)$ for a sufficiently
regular fixed function.  Note that a classical Kasner space-time
corresponds to the case of $\alpha \equiv 1$ (see \cite{KSS}
also).

By applying \emph{Lemma 2.8}, we obtained useful expressions for
the Ricci tensor and the scalar curvature of generalized
Robertson-Walker and generalized Kasner space-times. These
expressions allowed us to classify possible Einstein (respectively
with constant scalar curvature) generalized Kasner space-times of
dimension $4$. We also obtained some partial results for grater
dimensions.

Finally, in order to study curvature properties of multiply warped
product space-times associated to the BTZ $(2+1)$-dimensional
black hole solutions, we made applications of the previously
obtained curvature formulas. As a consequence, we characterized
the Einstein BTZ (respectively with constant scalar curvature), in
terms of the square lapse function.

\bigskip

In forthcoming papers we plan to focus on a specific
generalization of the structures studied here, which is
particularly useful in different fields such as relativity,
extra-dimension theories (Kaluza-Klein, Randall-Sundrum), string
and super-gravity theories, spectrum of Laplace-Beltrami operators
on $p$-forms, among others. Roughly speaking, we will consider a
mixed structure between a multiply warped product and a conformal
change in the base. Naturally, our main interest is the study of
curvature properties. As we have made progress on this subject, we
realized that these curvature related properties are interesting
and worth to study not only for the physical point of view (see
for instance, the several recent works of Gauntlett, Maldacena,
Argurio, Schmidt, among many others), but also for exclusive
nonlinear partial differential equations involved. Indeed, the
curvature related questions arise problems of existence,
uniqueness, bifurcation, study of critical points, etc. (see
Example 5.2 above and the different works of Aubin, Hebey, Yau,
Ambrosetti, Choquet-Bruat among others).

\bigskip

\begin{center}
{\bf Acknowledgements}
\end{center}

We would like to thank the referee for his valuable comments and
suggestions.

\end{document}